\documentclass{article}
\usepackage{amssymb}
\usepackage{graphicx}
\usepackage{amsmath}

\input{tcilatex}
\begin{document}

\begin{center}
\textbf{Free boundary value problems for abstract elliptic equations and
applications\ }\ 
\end{center}

\QTP{Body Math}
\ 

\begin{center}
{\textbf{Veli B. Shakhmurov}}

Department of Mechanical engineering, Okan University, Akfirat, Tuzla 34959
Istanbul, Turkey,

E-mail: veli.sahmurov@okan.edu.tr

\ \ \ \ \ \ \ \ \ \ \ \ \ \ 
\end{center}

\QTP{Body Math}
\ \ \ \ \ \ \ \ \ \ \ \ \ \ \ \ \ \ \ \ \ \ \ \ \ \ \ \ \ \ \ \ \ \ \ \ \ \
\ \ \ \ \ \ \ \ \ \ \ \ \ \ 

\begin{center}
\textbf{Abstract}
\end{center}

\begin{quote}
\ \ \ \ \ \ \ \ \ \ \ \ \ \ \ 
\end{quote}

Free bondary value problem for elliptic differential-operator equations with
variable coefficients is studied. The uniform maximal regularity properties
and Fredholmness of this problem are obtained in vector-valued Holder spaces.

\begin{center}
\qquad\ \ \ \ \ 

\end{center}


\begin{center}
\textbf{1. Introduction, notations and background }
\end{center}

In last years, the maximal regularity properties of boundary value problems
(BVPs) for differential-operator equations (DOEs) have found many
applications in PDE, psedo DE and in the different physical process (see for
references $\left[ 1-4\right] $, $\left[ 6\right] ,$ $\left[ 8\right] $, $%
\left[ 10\right] $, $\left[ 12-23\right] $, $\left[ 27-28\right] $ ).

Let $\Omega $ be a domin in $R^{n}$ and $E$ is a Banach space. $%
C_{b}^{\left( m\right) }\left( \Omega ;E\right) $\ will denote the spaces of 
$E$-valued bounded uniformly stongly continuous and $m$-times continuously
differentiable functions on $\Omega $. For $m=0$ it denotes by $C_{b}\left(
\Omega ;E\right) .$ Let $\mathbb{C}$ denote the set of complex numbers. For $%
E=\mathbb{C}$ the space $C^{\left( m\right) }\left( \Omega ;E\right) $ will
be denoted by $C_{b}^{\left( m\right) }\left( \Omega \right) .$ Moreover, $%
C_{b}^{\infty }\left( \Omega ;E\right) $ denotes spaces of $E$-valued
bounded strongly continiously differentiable functions of arbitrary order.
We put $\mathbb{R=}\left( -\infty ,\infty \right) $ and $\mathbb{R}%
_{+}=\left( 0,\infty \right) .$ Let $f\left( x\right) $ is a $E-$valued
function and $f\left( x\right) \neq 0$. Consider 
\begin{equation*}
\Omega _{f}=\left\{ \left( x,y\right) \in \mathbb{R}\times \mathbb{R}_{+},%
\text{ }f\in C_{b}\left( \mathbb{R};E\right) ,\text{ }0<y<\left\Vert f\left(
x\right) \right\Vert \right\} .
\end{equation*}

The boundaries of $\Omega _{f}$ are given by 
\begin{equation*}
\Gamma _{0}=\mathbb{R}\times 0,\text{ }\Gamma _{f}=\left\{ \left( x,y\right)
\in \mathbb{R}\times \mathbb{R}_{+},\text{ }y=\left\Vert f\left( x\right)
\right\Vert \right\} .
\end{equation*}%
Consider the following \ problem: Given $f_{0},$ $\nu \in C_{b}^{\left(
2\right) }\left( \mathbb{R};E\right) .$ Find a pair of functions $(u,f)$
possessing the regularity%
\begin{equation}
f\in C_{b}^{\left( 1\right) }\left( \left[ 0,\right. T;C_{b}^{\left(
1\right) }\left( \mathbb{R};E\right) \right) ,  \tag{1.1}
\end{equation}%
\begin{equation*}
u\left( t,.\right) \in W_{p}^{2}\left( \Omega _{f\left( t\right) };E\right) ,%
\text{ }t\in J=\left[ 0,\right. \left. T\right) .
\end{equation*}%
and satisfying the following equations a.e.

\begin{equation*}
-\Delta u\left( t,z\right) +A\left( x\right) u\left( t,z\right) =0\text{, }%
t\in J,\text{ }z\in \Omega _{f\left( t\right) }
\end{equation*}%
\begin{equation*}
\frac{\partial u}{\partial y}=0,\text{ }t\in J,\text{ }z\in \Omega _{f\left(
t\right) },
\end{equation*}%
\begin{equation}
u\left( t,z\right) =f\left( t,x\right) \text{, }t\in J,\text{ }z\in \Gamma
_{f\left( t\right) },  \tag{1.2}
\end{equation}%
\begin{equation*}
\lim\limits_{z\rightarrow \infty }u\left( t,z\right) =\nu \left( t\right) ,%
\text{ }t\in \left[ 0,\right. \left. T\right) ,
\end{equation*}%
\begin{equation*}
f_{t}\left( t,x\right) +\sqrt{1+f_{x}^{2}\left( t,x\right) }\frac{\partial }{%
\partial n}u\left( t,z\right) =0,\text{ }t\in \left( 0,T\right) ,\text{ }%
z\in \Gamma _{f\left( t\right) },
\end{equation*}%
\begin{equation*}
f\left( 0,x\right) =f_{0}\left( x\right) ,\text{ }x\in \mathbb{R},
\end{equation*}%
where $A$ is a linear operator in a Banach space $E$ and $z=\left(
x,y\right) $ represents a generic point in $\bar{\Omega}_{f}$. Moreover, $%
\Delta $ denotes the Laplace operator with respect to the Euclidean metric, $%
\frac{\partial }{\partial n}$ denotes the derivative in direction of the
outer unit normal $n$ at $\Gamma _{f\left( t\right) }$.

Maximal regularity properties of partial DOEs in $L_{p}$ spaces were studied
in $\left[ 1\right] $, $\left[ 4\right] $, $\left[ 7\right] $, $\left[ 18-23%
\right] .$ The results in $\left[ 4\right] $ and $\left[ 18-23\right] $ were
restricted to rectangular domain and equations that were not contained mixed
derivatives in leading part. Moreover, problems investigated in $\left[ 1%
\right] $ and $\left[ 8\right] $ involve only bounded operator coefficients.
In $\left[ 18\right] $ the Dirichlet\ problem for the elliptic
differential-operator equation of the second order in general domain was
studied.

In contrast to all above we study general BVP $\left( 1.1\right) $ for
equation with unbounded operator coefficients in the general domain.

Consider the BVP 
\begin{equation*}
Lu=\sum\limits_{i,j=1}^{n}a_{ij}\left( x\right) \frac{\partial ^{2}u}{%
\partial x_{i}\partial x_{j}}-A\left( x\right) u\left( x\right) =F\left(
x\right) ,
\end{equation*}%
\begin{equation}
L_{1}u=u\mid _{_{\Gamma }}=0,  \tag{1.3}
\end{equation}%
where $\Gamma $ is a boundary of region $G\subset R^{n}$ and $a_{ij}$ are
real-valued functions on $\bar{G}.$

We say that the problem $\left( 1.3\right) $ is maximal $H$-regular (or
separable in Holder space $C^{\gamma }$) if:

(1) for all $F\in C^{\gamma }\left( G;E\right) $ there exists a unique
solution $u\in $ $C^{2,\gamma }\left( G;E\left( A\right) ,E\right) $
satisfying $\left( 1.3\right) $ a.e. on $G;$

(2) there exists a positive constant $C$ independent of $F$ such that 
\begin{equation*}
\sum\limits_{i,j=1}^{n}\left\Vert \frac{\partial ^{2}u}{\partial
x_{i}\partial x_{j}}\right\Vert _{C^{\gamma }\left( G;E\right) }+\left\Vert
Au\right\Vert _{C^{\gamma }\left( G;E\right) }\leq C\left\Vert F\right\Vert
_{C^{\gamma }\left( G;E\right) }.
\end{equation*}%
Let $G$ denote the operator generated by the problem $\left( 1.3\right) $
for $\lambda =0$, i.e., 
\begin{equation*}
D\left( G\right) =C_{0}^{2,\gamma }\left( G;E\left( A\right) ,E\right)
=\left\{ u\in C^{2,\gamma }\left( G;E\right) \cap C\left( G;E\left( A\right)
\right) ,\right.
\end{equation*}%
\begin{equation*}
\left. u\mid _{_{\Gamma }}=0\right\} ,\text{ }Gu=Lu.\text{ }
\end{equation*}

The paper is organized as follows: Section1 collects definitions and
background materials, embedding theorems of Sobolev-Lions spaces.

Let $\mathbb{C}$ be the set of complex numbers and\ 
\begin{equation*}
S\left( \varphi \right) =\left\{ \text{\ }\lambda \in \mathbb{C}\text{, }%
\left\vert \arg \lambda \right\vert \leq \varphi \right\} \cup \left\{
0\right\} ,0\leq \varphi <\pi .
\end{equation*}

Let $E_{1}$ and $E_{2}$ be two Banach spaces. $L\left( E_{1},E_{2}\right) $
denotes the space of all bounded linear operators from $E_{1}$ to $E_{2}.$
For $E_{1}$ $=E_{2}=E$ it will be doneted by $L\left( E\right) .$

\ A linear operator\ $A$ is said to be $\varphi -$positive in a Banach\
space $E$ with bound $M>0$ if $D\left( A\right) $ is dense on $E$ and 
\begin{equation*}
\ \left\Vert \left( A+\lambda I\right) ^{-1}\right\Vert _{B\left( E\right)
}\leq M\left( 1+\left\vert \lambda \right\vert \right) ^{-1}
\end{equation*}%
with $\lambda \in S\left( \varphi \right) ,$ $\varphi \in \left[ 0,\right.
\left. \pi \right) $, where $I$ is an identity operator in $E$.

Sometimes instead of $A+\lambda I$\ will be written $A+\lambda $ and will
denoted by $A_{\lambda }.$ It is known that ($\left[ \text{25, \S 1.15.1}%
\right] )$ there exist fractional powers\ $A^{\theta }$ of positive operator 
$A.$ Let $E\left( A^{\theta }\right) $ denote the space $D\left( A^{\theta
}\right) $ with graphical norm 
\begin{equation*}
\left\Vert u\right\Vert _{E\left( A^{\theta }\right) }=\left( \left\Vert
u\right\Vert ^{p}+\left\Vert A^{\theta }u\right\Vert ^{p}\right) ^{\frac{1}{p%
}},1\leq p<\infty ,\text{ }-\infty <\theta <\infty .
\end{equation*}

A linear operator $A\left( x\right) $ is said to be positive in $E$
uniformly in $x$ if $D\left( A\left( x\right) \right) $ is independent of$\
x $, $D\left( A\left( x\right) \right) $ is dense in $E$ and 
\begin{equation*}
\left\Vert \left( A\left( x\right) +\lambda I\right) ^{-1}\right\Vert \leq
M\left( 1+\left\vert \lambda \right\vert \right) ^{-1}
\end{equation*}%
for all $\lambda \in S\left( \varphi \right) $ and $\varphi \in \left[
0,\right. \left. \pi \right) $.

Let $\Omega $ be a domain in $R^{n}$. $C(\Omega ,E)$ and $C^{m}(\Omega ;E)$
will denote the spaces of $E$-valued bounded uniformly strongly continuous
and $m$-times continuously differentiable functions on $\Omega ,$
respectively.

Let $0<\gamma \leq 1.$ $C^{\gamma }\left( \Omega ;E\right) $ denotes the
space of $E$-valued strongly bounded continiuous functions that are defined
on $\Omega \subset R^{n}$ with the norm

\begin{equation*}
\left\Vert f\right\Vert _{C^{\gamma }\left( \Omega ;E\right) }=\left\Vert
f\right\Vert _{C^{\gamma }\left( \Omega ;E\right) }+\left[ f\right] ^{\gamma
}\left( E\right) ,
\end{equation*}%
where 
\begin{equation*}
\left[ f\right] ^{\gamma }\left( E\right) =\sup\limits_{x\neq y,x,y\in
\Omega }\frac{\left\Vert f\left( x\right) -f\left( y\right) \right\Vert _{E}%
}{\left\vert x-y\right\vert ^{\gamma }}.
\end{equation*}

$C^{\gamma ,m}\left( \Omega ;E\right) $ denotes the space of $E$-valued
strongly bounded continiuous functions that are defined on $\Omega \subset
R^{n}$ with the norm

\begin{equation*}
C^{\gamma ,m}\left( \Omega ;E\right) =\left\{ f\in C^{m}\left( \Omega
;E\right) ,\text{ }f^{\left( m\right) }\in C^{\gamma }\left( \Omega
;E\right) \right\} ,
\end{equation*}

\begin{equation*}
\left\Vert f\right\Vert _{C^{\gamma ,m}\left( \Omega ;E\right) }=\left\Vert
f\right\Vert _{C^{m}\left( \Omega ;E\right) }+\left\Vert f^{\left( m\right)
}\right\Vert _{C^{\gamma }\left( \Omega ;E\right) }<\infty .
\end{equation*}

Let $E_{0}$ and $E$ be two Banach spaces and $E_{0}$ is continuously and
densely embedded into $E.$ Let $m$ be a natural number$.$

Let $\alpha =\left( \alpha _{1},\alpha _{2},...,\alpha _{n}\right) $ are $n$
tuples of nonnegative integer numbers and%
\begin{equation*}
D^{\alpha }=\frac{\partial ^{\left\vert \alpha \right\vert }}{\partial
x_{1}^{\alpha _{1}}\partial x_{2}^{\alpha _{2}}...\partial x_{n}^{\alpha
_{n}}}.
\end{equation*}%
\ $C^{m,\gamma }\left( \Omega ;E_{0},E\right) $ denote the space of $E_{0}$%
-valued bounded uniformly stongly continuous and $m$-times continuously
differentiable functions on $\Omega $ with norm%
\begin{equation*}
\left\Vert f\right\Vert _{C^{m,\gamma }\left( \Omega ;E_{0},E\right)
}=\left\Vert f\right\Vert _{C^{m,\gamma }\left( \Omega ;E\right)
}+\left\Vert f\right\Vert _{C^{\gamma }\left( \Omega ;E_{0}\right) }.
\end{equation*}

\ For $E_{0}=$ $E$ the space $C^{m,\gamma }\left( \Omega ;E_{0},E\right) $
will denoted by $C^{m,\gamma }\left( \Omega ;E\right) .$

Let $A$ be a linear operator in a Banach space $E$ so that is a generator of
analytic semigroup $U\left( t\right) =U_{A}\left( t\right) .$ Let%
\begin{equation*}
D_{A}\left( \theta ,p\right) =\left\{ u\in E,\text{ }\left\Vert u\right\Vert
_{D_{A}\left( \theta ,p\right) }=\left\Vert t^{1-\theta -\frac{1}{p}%
}AU\left( t\right) u\right\Vert _{L^{p}\left( 0,1;E\right) }<\infty \right\} 
\text{ }
\end{equation*}%
\begin{equation*}
\text{for }1\leq p<\infty \text{ and }\theta \in \left( 0,1\right) ;
\end{equation*}%
\begin{equation*}
D_{A}\left( \theta ,\infty \right) =\left\{ u\in E,\text{ }\left\Vert
u\right\Vert _{D_{A}\left( \theta ,p\right) }=\left\Vert t^{1-\theta
}AU\left( t\right) u\right\Vert _{L^{\infty }\left( 0,1;E\right) }<\infty
\right\} \text{ }
\end{equation*}

\begin{equation*}
\text{for }p=\infty \text{ and }0<\theta \leq 1;
\end{equation*}%
\begin{equation*}
D_{A}\left( \theta \right) =\left\{ u\in D_{A}\left( \theta ,\infty \right) ,%
\text{ }\lim\limits_{t\rightarrow 0}t^{1-\theta }AU\left( t\right)
u=0\right\} .\text{ }
\end{equation*}%
From $\left[ \text{15, Proposition 2.2.2}\right] $ we obtain the following%
\begin{equation*}
D_{A}\left( \theta ,p\right) =\left( E,D\left( A\right) \right) _{\theta ,p},%
\text{ for }1\leq p<\infty \text{, }\theta \in \left( 0,1\right) ,
\end{equation*}%
\begin{equation*}
D_{A}\left( \theta ,\infty \right) =\left( E,D\left( A\right) \right)
_{\theta ,\infty }\text{ for }0<\theta \leq 1,
\end{equation*}%
\begin{equation*}
D_{A}\left( \theta \right) =\left( E,D\left( A\right) \right) _{\theta },%
\text{ for }\theta \in \left( 0,1\right) .
\end{equation*}

$H\left( E_{0},E\right) $ denotes the class of linear operators that are
isomorphisim from $E_{0}$ onto $E$ and are negative generators of stronge
continious and analytic semigroups.

Let $S^{^{\prime }}\left( R^{n};E\right) $ denote the space of all\
continuous linear operators $L:S\left( R^{n};E\right) \rightarrow E,$
equipped with the bounded convergence topology. Recall $S\left(
R^{n};E\right) $ is norm dense in $L^{p}\left( R^{n};E\right) $ when $1\leq
p<\infty .$

Let $L_{q}^{\ast }\left( E\right) $ denote the space of all $E-$valued
functions $u\left( t\right) $ such that 
\begin{equation*}
\left\Vert u\right\Vert _{L_{q}^{\ast }\left( E\right) }=\left(
\int\limits_{0}^{\infty }\left\Vert u\left( t\right) \right\Vert _{E}^{q}%
\frac{dt}{t}\right) ^{\frac{1}{q}}<\infty \text{, }1\leq q<\infty
,\left\Vert u\right\Vert _{L_{\infty }^{\ast }\left( E\right)
}=\sup_{0<t<\infty }\left\Vert u\left( t\right) \right\Vert _{E}.
\end{equation*}%
Let $F$ denote the Fourier transform. Fourier-analytic definition of $E-$%
valued Besov space on $R^{n}$ are defined as in $\left[ \text{25 \S\ 3}%
\right] $, i.e.,%
\begin{equation*}
B_{p,q}^{s}\left( R^{n};E\right) =\left\{ u\in S^{^{\prime }}\left(
R^{n};E\right) ,\right. \text{ }
\end{equation*}%
\begin{equation*}
\left\Vert u\right\Vert _{B_{p,q}^{s}\left( R^{n};E\right) }=\left\Vert
F^{-1}\sum\limits_{k=1}^{n}t^{\varkappa _{k}-s_{k}}\left( 1+\left\vert \xi
_{k}\right\vert ^{\varkappa _{k}}\right) e^{-t\left\vert \xi \right\vert
^{2}}Fu\right\Vert _{L_{q}^{\ast }\left( L_{p}\left( R^{n};E\right) \right) }%
\text{,}
\end{equation*}%
\begin{equation*}
\left. p,\text{ }q\in \left[ 1,\infty \right] \text{, }\varkappa _{k}>s_{k},%
\text{ }s=\left( s_{1},s_{2},...,s_{n}\right) \right\} .
\end{equation*}

For appropriate domain $\Omega \subset R^{n}$ the space $B_{p,q}^{s}\left(
\Omega ;E\right) $ is defined as usulal restriction of the space $%
B_{p,q}^{s}\left( R^{n};E\right) .$

For $E=\mathbb{C}$ the space $B_{p,q}^{s}\left( \Omega ;\mathbb{C}\right) $
will be denoted by $B_{p,q}^{s}\left( \Omega \right) $.

Let $h^{s}=h^{s}\left( R^{n};E\right) $ denote the closure of $S\left(
R^{n};E\right) $ in $B_{\infty ,\infty }^{s}\left( R^{n};E\right) $. Assume
that $\Omega $ is an open subset of $R^{n}$ and let $r_{\Omega }$ denote the
restriction operator with respect to $\Omega $, i.e., $r_{\Omega }u=u\mid
_{\Omega }$for $B_{\infty ,\infty }^{s}\left( R^{n};E\right) .$ Here, $%
h^{s}\left( \Omega ;E\right) $ is defined as the closure of $r_{\Omega
}\left( S\left( R^{n};E\right) \right) $ in $B_{\infty ,\infty }^{\gamma
}\left( \Omega ;E\right) $ and $h^{m,\gamma }\left( \Omega ;E\right) $ is
defined as the closure of $r_{\Omega }\left( S\left( \Omega ;E\right)
\right) $ in $C^{m,\gamma }\left( \Omega ;E\right) .$ For $E=\mathbb{C}$ the
spaces $h^{s}\left( \Omega ;E\right) ,$ $h^{m,\gamma }\left( \Omega
;E\right) $ will be denoted by $h^{s}$ and $h^{m,\gamma }$, respectively.
Moreover, let $C_{0}^{s}\left( \Omega ;E\right) $ denote the closure of $%
C^{\infty }\left( \Omega ;E\right) $ in $C^{s}\left( \Omega ;E\right) .$

\bigskip Let 
\begin{equation*}
h_{+}^{s}=h_{+}^{s}\left( \mathbb{R};E\right) =\left\{ g\in h^{s}\text{, }%
g\left( x\right) \neq 0\right\} ,
\end{equation*}%
\begin{equation*}
X=h^{\alpha }\left( Q;E\right) \text{, }X_{k}=h^{k,\alpha }\left( Q;E\right) 
\text{, }Y=h^{2,\alpha }\left( Q;E\left( A\right) ,E\right) =
\end{equation*}

\begin{equation*}
\left\{ u\in h^{\alpha }\left( Q;D\left( A\right) \right) \cap h^{2,\alpha
}\left( Q;E\right) \right\}
\end{equation*}%
with the norm%
\begin{equation*}
\left\Vert u\right\Vert _{Y}=\left\Vert u\right\Vert _{h^{\alpha }\left(
Q;D\left( A\right) \right) }+\left\Vert u^{\left( 2\right) }\right\Vert
_{h^{\alpha }\left( Q;E\right) }<\infty .
\end{equation*}

Here, 
\begin{equation*}
B_{+}^{2+\alpha }=\left\{ u\in B_{\infty ,\infty }^{2+\alpha }\left(
R^{n};E\right) \text{, }g\left( x\right) \neq 0\right\} \text{ }
\end{equation*}%
and 
\begin{equation*}
h^{m,\alpha }\left( A\right) =h^{m,\alpha }\left( \mathbb{R};D_{A}\left(
\alpha ,\infty \right) \right) \text{, }\alpha \in \left( 0,1\right) .
\end{equation*}

\textbf{Remark 1.1. }In order to formulate our result, let 
\begin{equation*}
h_{\nu }^{s}=h_{\nu }^{s}\left( E\right) =\left\{ \nu +g;\text{ }g\in
h^{s}\left( \mathbb{R};E\right) \right\} ,\text{ }h_{\nu }^{s,\alpha
}=\left\{ \nu +g;\text{ }g\in h^{s,\alpha }\left( \mathbb{R};E\right)
\right\}
\end{equation*}%
and $f\in h_{\nu }^{2,\alpha }$ given. Let $u_{f}$ denote the unique
solution of the BVP 
\begin{equation*}
-\Delta u+Au=0,\text{ }\partial _{y}u=0\text{ on }\Gamma _{0},\text{ }u=f%
\text{ on }\Gamma _{f},
\end{equation*}%
where $A$ is a linear operator in a Banach space $E.$%
\begin{equation*}
k_{f}=\frac{\left\Vert f\right\Vert ^{2}}{\left( 1+\left\Vert f\right\Vert
+\left\Vert f_{x}\right\Vert ^{2}\right) \left( 1+\left\Vert
f_{x}\right\Vert ^{2}\right) }
\end{equation*}%
and define \ 
\begin{equation*}
V_{\nu }=\left\{ f\in C_{b}^{2}\left( \mathbb{R};E\right) \right\} \text{, }%
f\left( x\right) \neq 0,\text{ }\partial _{y}u_{f}\left( x,f\left( x\right)
\right) <k_{f}\text{ for }x\in \mathbb{R}.
\end{equation*}

It is clear to see that $u_{\nu }\equiv \nu .$ Hence, $\nu \in V_{\nu }.$
More precisely, by following $\left[ \text{9, Lemma 5.10}\right] $ it can be
shown that $V_{\nu }$ is a open neighborhood of $\nu $ in $h_{\nu
}^{s,2}\left( E\right) $ and that%
\begin{equation*}
\text{diam}_{h^{2,s}}\left( V_{\nu }\right) =\sup\limits_{g,\upsilon \in
V_{\nu }}\left\Vert g-\upsilon \right\Vert _{h^{s,2}}=\infty .
\end{equation*}

Suppose now that$(u,f)$ is a classical solution of $\left( 1.1\right)
-\left( 1.2\right) $. We call $(u,f)$ a classical H\"{o}lder solution on $J$
if it possesses the additional regularity 
\begin{equation*}
f\in C\left( J;V_{\nu }\right) \cap C^{1}\left( J;h_{\nu }^{1,\alpha
}\right) \text{, }u\left( t,.\right) \in h_{\nu }^{2,\alpha }\left( \mathbb{R%
};E\right) \text{, }t\in J.
\end{equation*}

We will prove the following main result

\textbf{Theorem 1. }Given $f_{0}\in V_{\nu }$, there exist $%
t^{+}=t^{+}\left( f_{0}\right) $ and a unique maximal classical H\"{o}lder
solution$(u,f)$ of problem $\left( 1.1\right) -\left( 1.2\right) $ on $%
[0,t^{+})$. Moreover, the mapping $(t,f_{0})\rightarrow f$ defines a local $%
C^{\infty }$-semiflow on $V_{\nu }$. If $t^{+}<\infty $ and $f:$ $\left[
\left. 0,t^{+}\right) \right. \rightarrow V_{\nu }$ is uniformly continuous
then either 
\begin{equation*}
\lim\limits_{t\rightarrow t^{+}}\left\Vert f\left( t,.\right) \right\Vert
_{h^{2,\alpha }}=\infty \text{, }\lim\limits_{t\rightarrow
t^{+}}\inf\limits_{\upsilon \in V_{\nu }}\left\Vert f\left( t,.\right)
-\upsilon \right\Vert _{h^{2,\alpha }}=0\text{.}
\end{equation*}

In the first stage, we transform problem $\left( 1.1\right) -\left(
1.2\right) $ into a nonlinear problem on a fixed domain 
\begin{equation*}
\frac{df}{dt}+O\left( f\right) =0,\text{ }f\left( 0\right) =f_{0}
\end{equation*}%
with respect to only the unknown function $f$, which determines the free
boundary $\Gamma _{f},$ where $O$ is a nonlinear operator in $E.$

Then, by using the solution of the above problem we will show the exsistence
of regular solution of the free BVP $\left( 1.1\right) -\left( 1.2\right) .$

\begin{center}
\ \textbf{2. Transformed problem}
\end{center}

\bigskip Let $\nu =\nu \left( t\right) >0$ be fixed. Define 
\begin{equation*}
G_{\nu }=\left\{ g\in C_{b}^{\left( 2\right) }\left( \mathbb{R};E\right) ,%
\text{ }\nu \left( t\right) I+g\left( x\right) \neq 0\right\} ,
\end{equation*}%
where $I$ is an identity element in the Banach space $E.$

Consider the following transformation%
\begin{equation}
\left( x,y\right) =\varphi \left( x^{\prime },y^{\prime }\right) =\varphi
_{g}\left( x^{\prime },y^{\prime }\right) =\left( x^{\prime },1-\frac{%
y^{\prime }}{\nu +g\left( x^{\prime }\right) }\right) \text{, for }\left(
x^{\prime },y^{\prime }\right) \in \Omega _{f}.  \tag{2.1}
\end{equation}

It is easily verified that $\varphi _{g}$ is a diffeomorphism of class $%
C^{2} $ which maps $\Omega _{g}$ onto the strip $Q=\mathbb{R}\times \left(
0,1\right) $. Moreover,%
\begin{equation*}
\left( x^{\prime },y^{\prime }\right) =\varphi ^{-1}\left( x,y\right)
=\varphi _{g}^{-1}\left( x,y\right) =\left( x,\left( 1-y\right) g\left(
x\right) \right) \text{ for }\left( x,y\right) \in Q.
\end{equation*}

\begin{equation}
\varphi _{\ast }u=\varphi _{\ast }^{g}u=u\left( \varphi _{g}^{-1}\left(
x,y\right) \right) \text{ for }u\in W_{p}^{2}\left( \Omega _{g};E\left(
A\right) ,E\right) .\text{ }  \tag{2.2}
\end{equation}%
Let $u$ be an $E$-valued function defined on $Q$. Here, $u_{\mid _{\Gamma
_{i}}}$denote the restriction of $u$ on $\Gamma _{i}$, where 
\begin{equation*}
\Gamma _{i}=\mathbb{R}\times \left\{ i\right\} ,\text{ }\gamma _{i}u=u_{\mid
_{\Gamma }},\text{ }i=0,1.
\end{equation*}%
\textbf{Lemma 2.1. }Given $g\in \Phi _{\eta }$ and $\upsilon \in C^{2,\gamma
}\left( Q;E\left( A\right) ,E\right) ;$ under the the maps $\left(
2.2\right) $ the operators in $\left( 1.1\right) $ are transformed into the
following:%
\begin{equation}
B\left( g\right) \upsilon =-\varphi ^{g}\ast \left( \Delta +A\right) \left(
\varphi _{g}\ast \upsilon \right) \text{ on }\left[ 0,\right. \left.
T\right) \times Q,  \tag{2.3}
\end{equation}%
\begin{equation*}
B_{i}\left( g\right) \upsilon =\varphi _{\ast }^{g}\left( \left( \nabla
\left( \varphi _{g}^{\ast }\right) ,n_{i}\right) \right) \text{ on }\left(
0,T\right) \times \Gamma _{i},
\end{equation*}%
where $n_{0}=\left( -g_{x},1\right) $, $n_{1}=\left( 0,-1\right) $ denote
the outer normals according to $\Gamma _{f}$ and $\Gamma _{0}$, i.e.,

\begin{equation}
B_{0}\left( g\right) \upsilon =\varphi _{\ast }^{g}\left( -\frac{\partial }{%
\partial x}\varphi _{g}^{\ast }\upsilon +\frac{\partial }{\partial y}\varphi
_{g}^{\ast }\upsilon \right) \text{ on }\left( 0,T\right) \times \Gamma _{0},
\tag{2.4}
\end{equation}%
\begin{equation*}
B_{1}\left( f\right) \upsilon =\varphi _{\ast }^{g}\left( -\frac{\partial }{%
\partial y}\varphi _{g}^{\ast }\upsilon \right) \text{ on }\left[ 0,\right.
\left. T\right) \times \Gamma _{1}.
\end{equation*}

\textbf{Lemma 2.2. }Given $g\in \Phi _{\nu }$ and $\upsilon \in C^{2,\gamma
}\left( Q;E\left( A\right) ,E\right) .$ Under the map $\left( 2.2\right) $
the problem $\left( 1.1\right) $ is transformed into the following:%
\begin{equation*}
B\left( f\right) \upsilon =-\varphi ^{f}\ast \left( \Delta +A\right) \left(
\varphi _{g}\ast \upsilon \right) =0\text{ on }\left[ 0,\right. \left.
T\right) \times Q,
\end{equation*}%
\begin{equation*}
\upsilon =f\text{ \ on }\left[ 0,\right. \left. T\right) \times \Gamma _{0},
\end{equation*}%
\begin{equation}
B_{1}\left( f\right) \upsilon =0\text{ on }\left[ 0,\right. \left. T\right)
\times \Gamma _{1},  \tag{2.5}
\end{equation}%
\begin{equation*}
\lim\limits_{z\rightarrow \infty }\upsilon \left( t,z\right) =0,\text{ on }%
\left[ 0,\right. \left. T\right) ,
\end{equation*}%
\begin{equation*}
\frac{\partial g}{\partial t}+B_{0}\left( g\right) \upsilon =0\text{ on }%
\left( 0,T\right) \times \Gamma _{0},
\end{equation*}%
\begin{equation*}
g\left( 0,.\right) =g_{0}\left( x\right) \text{ on }\mathbb{R},
\end{equation*}%
\textbf{\ }A pair $\left( \upsilon ,g\right) $ is called a solution of the
problem $\left( 2.5\right) $ if%
\begin{equation}
g\in C_{b}\left( \left[ 0,\right. \left. T\right) ;\Phi \right) \cap
C_{b}^{\left( 1\right) }\left( \left[ 0,\right. \left. T\right)
;C_{b}^{\left( 1\right) }\left( \mathbb{R}\right) \right) ,  \tag{2.5}
\end{equation}%
\begin{equation*}
u\left( t,.\right) \in W_{p}^{2}\left( \Omega _{f\left( t\right) };E\right) ,%
\text{ }t\in \left[ 0,\right. \left. T\right)
\end{equation*}%
and $\left( \upsilon ,g\right) $ satisfies $\left( 2.5\right) $ a.e. on $%
\left[ 0,\right. \left. T\right) \times Q.$

\textbf{Condition 2.1. }Assume the following conditions are satisfied:

(1) $\dsum\limits_{i,j=1}^{2}a_{ij}\xi _{i}\xi _{j}\geq C\left\vert \xi
\right\vert ^{2}$, for $\xi =\left( \xi _{1},\xi _{2},...,\xi _{n}\right)
\in R^{n}$ and $C>0;$

(2) operator $A$ is a positive operator in a Banach space $E$ for some $%
\varphi \in \left( 0\right. ,\left. \pi \right] .$

In a similar way as in $\left[ \text{9, Lemma 2.2}\right] $ we obtain

\textbf{Lemma 2.3. }Assume the Condition 2.1 are satisfied. Then for given $%
g\in \Phi ,$ we have%
\begin{equation}
B\left( g\right) u=\dsum\limits_{i,j=1}^{2}-a_{ij}\left( g\right) \frac{%
\partial ^{2}u}{\partial x_{i}\partial x_{j}}+a_{2}\left( g\right) \frac{%
\partial u}{\partial x_{2}}+A_{g}u,  \tag{2.6}
\end{equation}%
\begin{equation*}
B_{i}\left( g\right) u=\dsum\limits_{j=1}^{2}b_{ji}\left( g\right) \gamma
_{i}\frac{\partial u}{\partial x_{j}}\text{, }i=0,1,
\end{equation*}%
and 
\begin{equation*}
\dsum\limits_{i,j=1}^{2}a_{ij}\left( g\right) \xi _{i}\xi _{j}\geq \alpha
\left( g\right) \left\vert \xi \right\vert ^{2}\text{, for }\xi =\left( \xi
_{1},\xi _{2},...,\xi _{n}\right) \in R^{n},
\end{equation*}%
where $\alpha \left( g\right) >0$, $\gamma _{i}$ are trace operators from $Q$
to $\Gamma _{i}$, $i=0,1,$ 
\begin{equation*}
a_{11}\left( g\right) =1\text{, }a_{12}\left( g\right) =a_{21}\left(
g\right) =\frac{\beta g_{x_{1}}}{\nu +g},\text{ }a_{22}\left( g\right) =%
\frac{1+\beta ^{2}g_{x_{1}}^{2}}{\left( \nu +g\right) ^{2}},\text{ }\beta
=1-x_{1},
\end{equation*}%
\begin{equation}
a_{2}\left( g\right) =\frac{\beta }{\nu +g}\left[ \frac{2g_{x_{1}}^{2}}{\nu
+g}-g_{x_{1}x_{1}}\right] ,\text{ }b_{10}\left( g\right) =-g_{x_{1}},\text{ }%
b_{20}\left( g\right) =-\frac{1+g_{x_{1}}^{2}}{\nu +g},  \tag{2.7}
\end{equation}%
\begin{equation*}
b_{11}\left( g\right) =0\text{, }b_{21}\left( g\right) =\frac{1}{\nu +g}%
\text{, }\alpha \left( g\right) =\frac{1}{1+\left( \nu +g\right) ^{2}+\beta
^{2}g_{x_{1}}^{2}},\text{ }A\left( g\right) =A\left( \varphi _{g}\right) .
\end{equation*}

\begin{center}
\textbf{3. Abstract elliptic equation in the fixed domain}
\end{center}

\bigskip In this section we study the elliptic BVP 
\begin{equation}
B\left( g\right) u=-\dsum\limits_{i,j=1}^{2}a_{ij}\left( g\right) \frac{%
\partial ^{2}u}{\partial x_{i}\partial x_{j}}+a_{2}\left( g\right) \frac{%
\partial u}{\partial x_{2}}+A\left( x\right) u=f,  \tag{3.1}
\end{equation}%
\begin{equation}
B_{i}\left( g\right) u=\dsum\limits_{j=1}^{2}b_{ji}\left( g\right) \gamma
_{i}\frac{\partial u}{\partial x_{j}}=f_{i}\text{, }i=0,1,  \tag{3.2}
\end{equation}%
where $B\left( g\right) $ and $B_{i}\left( g\right) $ are differential
operators defined by $\left( 2.6\right) .$

We will derive a priori estimates as well as isomorphism properties in
framwork of abstract Holder spaces.

\textbf{Condition 3.1. }Assume the following conditions are satisfied:

(1) $a_{ij}\in C^{0,\alpha }\left( \bar{G}\right) ,$ $a_{ij}=a_{ji};$

(2) $\dsum\limits_{i,j=1}^{2}a_{ij}\left( g\right) \xi _{i}\xi _{j}\geq
C\left( g\right) \left\vert \xi \right\vert ^{2}$, for $\xi =\left( \xi
_{1},\xi _{2},...,\xi _{n}\right) \in R^{n}$ and $C\left( g\right) >0;$

(3) operator $A\left( g\right) $ is uniformly positive in a Banach algebra $%
E $ for some $\varphi \in \left( 0\right. ,\left. \pi \right] .$

Here $\partial B\left( g\right) \left[ \psi ,.\upsilon \right] $ denotes the
Gateaux derivative of operator function $B\left( g\right) $ at $\psi $ in
the direction of $\upsilon .$

\textbf{Lemma 3.1. }Suppose the Condition 3.1 is satisfied and $A\left(
g\right) $ is Gateaux differentiable for $g\in h_{\Phi }^{2,\alpha }$, $%
\alpha \in \left( 0,1\right) $. Then the map $g\rightarrow O_{0}\left(
g\right) =\left\{ B\left( g\right) ,\text{ }B_{0}\left( g\right) \right\} $
is bounded linear operator-function from $Y$ into $X\times h^{1,\alpha
}\left( A\right) $ and have continious derivatives of all order with respect
to $g\in h_{\Phi }^{2,\alpha }$, i.e. 
\begin{equation*}
B\left( .\right) \in C^{\infty }\left( h_{\Phi }^{2,\alpha };L\left(
Y,X\right) \right) ,\text{ }B_{0}\left( .\right) \in C^{\infty }\left(
h_{\Phi }^{2,\alpha };L\left( Y,E\right) \right)
\end{equation*}%
and 
\begin{equation*}
\partial B\left( g\right) \left[ u,\upsilon \right] =\left( \partial B\left(
g\right) u\right) \upsilon =\frac{2\beta }{\nu +g}\left\{ \left( \frac{g_{x}u%
}{\nu +g}-u_{x}\right) \upsilon _{x_{1}x_{2}}+\frac{u}{\left( \nu +g\right)
^{2}}\left( \left( \frac{1}{\beta }+\beta g_{x}^{2}\right) \right. \right. -
\end{equation*}%
\begin{equation*}
\left. \frac{\beta }{\nu +g}g_{x}u_{x}\right) \upsilon _{x_{2}x_{2}}-\left.
\left( \frac{g_{x}u}{\left( \nu +g\right) ^{2}}-\frac{g_{xx}u+4g_{x}u_{x}}{%
2\left( \nu +g\right) }+\frac{u_{xx}}{2}\right) \upsilon _{x_{2}}\right\}
+\partial A\left( g\right) \left[ u,\upsilon \right] ,
\end{equation*}%
\begin{equation*}
\partial B_{0}\left( g\right) \left[ u,\upsilon \right] =-u_{x}\upsilon
_{x_{1}}+\frac{\beta }{\nu +g}\left[ \frac{u\left( 1+g_{x}^{2}\right) }{\nu
+g}-2g_{x}u_{x}\right] \upsilon _{x_{2}}
\end{equation*}%
for $g\in h_{+}^{2,\alpha }$and $u,$ $\upsilon \in Y.$

\textbf{Proof. }$\ $It is clear to see that $\left( \varphi ,\upsilon
\right) \rightarrow \varphi \upsilon $ is bilinear and continuous from $%
h^{s} $ into $Y.$ Moreover, the mapping 
\begin{equation*}
g\rightarrow \frac{1}{\nu +g}u
\end{equation*}
are continious and are infinitely many differentible from $h^{s}$ into $Y.$
By using the definition of the space $Y$ and Lemma 2.3 we get that for all
fixed $g\in h_{\Phi }^{2,\alpha }$ the operator $u\rightarrow B\left(
g\right) $ is bounded linear operator from $Y$ into $X.$ So, we obtain that 
\begin{equation*}
B\left( .\right) \in C^{\infty }\left( h_{\Phi }^{2,\alpha },L\left(
Y,X\right) \right) .
\end{equation*}

Hence, in view of Lemma 2.3 we obtain 
\begin{equation*}
B_{0}\left( .\right) \in C^{\infty }\left( h_{\Phi }^{2,\alpha },L\left(
Y,E\right) \right) .
\end{equation*}

By using $\left[ \text{4, Theorem 2}\right] $ we obtain the following:

\textbf{Theorem 3.1. }Suppose the Condition 3.1 is satisfied and $\alpha \in
\left( 0,1\right) $. Then for $\lambda \in S\left( \varphi \right) $ and for
sufficiently large $\left\vert \lambda \right\vert :$

(a) the operator $u\rightarrow \tilde{O}\left( g\right) u=\left\{ \left(
B\left( g\right) +\lambda \right) u,\text{ }\gamma _{0}u,\text{ }\left( \eta
+g\right) B_{1}\left( g\right) \right\} $ is isomorphism from $Y$ onto $%
X\times h^{2,\alpha }\left( A\right) \times h^{1,\alpha }\left( A\right) ;$

(b) for $\mu >0$ the operator $u\rightarrow \left\{ \left( B\left( g\right)
+\lambda \right) u,\text{ }\mu \gamma _{0}u,\text{ }\left( \eta +g\right)
B_{1}\left( g\right) \right\} $ is isomorphism from $Y$ onto $X\times
h^{2,\alpha }\left( A\right) \times h^{1,\alpha }\left( A\right) ;$

(c) For $u\in Y$ there exists a positive constant $C$, depending only on $g$%
, $\eta $, $p$ and $E$ such that the coercive estimate holds:%
\begin{equation}
\left\Vert u\right\Vert _{Y}\leq C\left( \left\Vert \left( B\left( g\right)
+\lambda \right) u\right\Vert _{X}+\left\Vert \gamma _{0}u\right\Vert
_{h^{2,\alpha }\left( E\right) }+\left\Vert \left( \nu +g\right)
B_{1}u\right\Vert _{h^{1,\alpha }\left( E\right) }\right) .  \tag{3.3}
\end{equation}

\textbf{Proof}. Indeed, since the domain $Q$ is a strip, the functios $g$, $%
\eta $ are fixed smooth functions, by virtue of trace theorem in Hoder space 
$\left[ \text{15, \S\ 2}\right] $ we obtain the assertion.

Consider now, the following BVPs%
\begin{equation}
-\dsum\limits_{i,j=1}^{2}a_{ij}\left( g\right) \frac{\partial ^{2}u}{%
\partial x_{i}\partial x_{j}}+a_{2}\left( g\right) \frac{\partial u}{%
\partial x_{2}}+A\left( g\right) u=F,  \tag{3.4}
\end{equation}%
\begin{equation}
B_{i}u=\dsum\limits_{j=1}^{2}b_{ji}\left( g\right) \gamma _{i}\frac{\partial
u}{\partial x_{j}}=0\text{, }i=0,1,  \tag{3.5}
\end{equation}%
\begin{equation}
-\dsum\limits_{i,j=1}^{2}a_{ij}\left( g\right) \frac{\partial ^{2}u}{%
\partial x_{i}\partial x_{j}}+a_{2}\left( g\right) \frac{\partial u}{%
\partial x_{2}}+A\left( g\right) u=0,  \tag{3.6}
\end{equation}%
\begin{equation}
\text{ }\dsum\limits_{j=1}^{2}b_{j0}\left( g\right) \gamma _{0}\frac{%
\partial u}{\partial x_{j}}=\psi \text{, }\dsum\limits_{j=1}^{2}b_{j1}\left(
g\right) \gamma _{1}\frac{\partial u}{\partial x_{j}}=0\text{.}  \tag{3.7}
\end{equation}%
Consider the operators $\tilde{S}\left( g\right) $ and $\tilde{K}\left(
g\right) $ generated by problems $\left( 3.4\right) -\left( 3.5\right) $ and 
$\left( 3.6\right) -\left( 3.7\right) $, respectively, i.e. 
\begin{equation*}
D\left( \tilde{S}\left( g\right) \right) =\left\{ u\in Y,\text{ }B_{i}u=0%
\text{, }i=0,1\right\} ,\text{ }
\end{equation*}

\begin{equation*}
\tilde{S}\left( g\right) u=-\dsum\limits_{i,j=1}^{2}a_{ij}\left( g\right) 
\frac{\partial ^{2}u}{\partial x_{i}\partial x_{j}}+a_{2}\left( g\right) 
\frac{\partial u}{\partial x_{2}}+A\left( g\right) u
\end{equation*}%
and 
\begin{equation*}
D\left( \tilde{K}\left( g\right) \right) =\left\{ u\in Y,\text{ }Bu=0,\text{ 
}B_{1}u=0\right\} \text{, }\tilde{K}\left( g\right) u=B_{0}u\in B_{0p}.
\end{equation*}

From the Theorem 3.5 we obtain that the inverse operators $\tilde{O}%
^{-1}\left( g\right) ,$ $\tilde{S}^{-1}\left( g\right) ,$ $\tilde{K}%
^{-1}\left( g\right) $ are bounded from $X\times h^{2,\alpha }\left(
A\right) \times h^{1,\alpha }\left( A\right) $, $X,$ $h^{2,\alpha }\left(
E\right) $ into $Y,$ respectively.

Here, 
\begin{equation*}
O\left( g\right) =\tilde{O}^{-1}\left( g\right) ,\text{ }S\left( g\right) =%
\tilde{S}^{-1}\left( g\right) \text{, }K\left( g\right) =\tilde{K}%
^{-1}\left( g\right) .
\end{equation*}

Assume $g\in h_{+}^{2,\alpha }\left( \mathbb{R};E\right) ,$ $\psi \in
D_{A}\left( \alpha \right) $, and put $u=$ $K\left( g\right) \upsilon .$
Then $u$ is the solution of the BVP $\left( 3.6\right) -\left( 3.7\right) .$

\textbf{Condition 3.1. }Assume $A\left( g\right) $ is Gateaux differentiable
for $g\in h_{+}^{2,\alpha }$, $\alpha \in \left( 0,1\right) $ and the
operator $\partial A\left( g\right) $ is uniformly $R$-positive in UMD (see
e.g $\left[ 8\right] $ for definitions) Banach algebra $E.$

\bigskip \textbf{Lemma 3.2. }Suppose the conditions 3.1 and 3.2 are
satisfied. Then we have 
\begin{equation*}
K\left( .\right) \in C^{\infty }\left( h_{+}^{2+\alpha };L\left( E,Y\right)
\right)
\end{equation*}%
and 
\begin{equation*}
\partial K\left( g\right) \left[ \upsilon ,\psi \right] =-S\left( g\right)
\partial B\left( g\right) \left[ \upsilon ,K\left( g\right) \psi \right]
\end{equation*}%
for $g\in h_{+}^{2,\alpha }$, $\psi \in E$ and $\upsilon \in Y.$

\textbf{Proof. }It follows from Lemma 3.1 and Theorem 3.1 that the map $%
g\rightarrow \tilde{O}\left( g\right) $ is isomorphism from $Y$ onto $%
X\times h^{2,\alpha }\left( A\right) \times h^{1,\alpha }\left( A\right) $
and have continious derivatives of all order with respect to $g\in
h_{+}^{2,\alpha }$, i.e. 
\begin{equation*}
\tilde{O}\in C^{\infty }\left( h_{+}^{2,\alpha };\text{Isom}Y,X\times
h^{2,\alpha }\left( A\right) \times h^{1,\alpha }\left( A\right) \right)
\end{equation*}%
and 
\begin{equation*}
\partial \tilde{O}\left( g\right) \psi =\left\{ \partial B\left( g\right) %
\left[ \psi ,.\right] ,\text{ }0\text{, }0\right\} \text{ for }\psi \in
h^{2,\alpha }\left( A\right) .
\end{equation*}

Then by reasoning as the Lemma 3.4 in $\left[ 6\right] $ we obtain the
assertion.\bigskip

\begin{center}
\textbf{4. The nonlinear operator for free BVPs}
\end{center}

\bigskip In this section we introduce the basic nonlinear operator and we
derive some properties of it.

Moreover, we show that the corresponding evolution problem involving this
operator is equivalent to the original problem $\left( 1.2\right) .$ Given $%
g\in h_{+}^{2,\alpha },$ we define the following operator%
\begin{equation*}
O\left( g\right) =B_{0}\left( g\right) K\left( g\right) g.
\end{equation*}

From lemmas $3.1$ and $3.2$ we get that 
\begin{equation}
O\in C^{\infty }\left( h_{+}^{2,\alpha },h^{1,\alpha }\left( E\right)
\right) .  \tag{4.1}
\end{equation}

Assume that $g_{0}\in h_{+}^{2,\alpha }$ , and let $\sigma =\left[ 0,\right.
\left. T\right) .$ A function $\sigma \rightarrow B_{10}$ is a classical
solution of 
\begin{equation}
\frac{dg}{dt}+O\left( g\right) =0,\text{ }g\left( 0\right) =g_{0}  \tag{4.2}
\end{equation}

iff $g\in C\left( \sigma ;h_{+}^{2,\alpha }\right) \cap $ $C^{1}\left(
\sigma ;h^{2,\alpha }\left( E\right) \right) $ and $g$ satisfies $\left(
4.2\right) $ pointwise.

\textbf{Lemma 4.1. }Suppose the Condition 3.1 is satisfied. Then for $%
g_{0}\in h_{+}^{2,\alpha }$:

(a) Suppose that $g$ is a classical solution of problem $\left( 4.2\right) $
on $\sigma $ and let $\upsilon \left( t,.\right) =K\left( g\left( t,.\right)
\right) g\left( t,.\right) .$ Then the pair $(\upsilon ,$ $g)$ is a
classical solution of $\left( 2.5\right) $ on $\sigma $, having the
additional regularity\ 
\begin{equation*}
g\in C\left( \sigma ;h_{U}^{2,\alpha }\right) \cap C^{1}\left( \sigma
;h^{1,\alpha }\right) ,
\end{equation*}%
\begin{equation}
\upsilon \left( t,.\right) \in Y\text{, }t\in \sigma ;  \tag{4.3}
\end{equation}

(b) Suppose that $(\upsilon ,$ $g)$ is a classical solution of $\left(
2.5\right) $ on $\sigma $, having the regularity $\left( 4.3\right) .$ Then $%
g$ is a classical solution of $\left( 4.2\right) $ on $\sigma $.

\textbf{Proof. }The proof is obtained from Lemma 2.2 and definitions of the
spaces $B_{ip},$ $i=0,1$ and $Y$.

For fixed $g\in h_{\Phi }$ consider the operator 
\begin{equation}
\upsilon \rightarrow B_{0}\left( g\right) K\left( g\right) \upsilon . 
\tag{4.4}
\end{equation}

In view of Theorem 3.1 we obtain that%
\begin{equation}
B_{0}\left( g\right) K\left( g\right) \in L\left( h^{2,\alpha }\left(
A\right) ,h^{1,\alpha }\left( A\right) \right) .  \tag{4.5}
\end{equation}

\textbf{Lemma 4.2. }Suppose the conditions 3.1 and 3.2 are satisfied. Then $%
O\in C^{\infty }\left( h_{+}^{2,\alpha },h^{1,\alpha }\left( A\right)
\right) $ and 
\begin{equation}
\partial O\left( g\right) \psi =B_{0}\left( g\right) K\left( g\right) \psi
+\partial B_{0}\left( g\right) \left[ \upsilon ,K\left( g\right) \psi \right]
-  \tag{4.6}
\end{equation}

\begin{equation*}
B_{0}\left( g\right) S\left( g\right) \partial B\left( g\right) \left[
\upsilon ,K\left( g\right) \psi \right]
\end{equation*}%
for $g\in h_{+}^{2,\alpha }$, $\psi \in h^{2,\alpha }\left( A\right) $ and $%
\upsilon \in Y.$

\textbf{Proof. }The assertion is obtained from Lemma 3.1 and Lemma 3.2.

\begin{center}
\bigskip \textbf{5. Linear equation with constant coefficients}
\end{center}

We put 
\begin{equation*}
R_{+}^{2}=\left\{ x=\left( x_{1},x_{2}\right) \in R^{2},\text{ }%
x_{2}>0\right\} .
\end{equation*}%
Consider the problem 
\begin{equation}
-\dsum\limits_{i,j=1}^{2}a_{ij}\frac{\partial ^{2}u}{\partial x_{i}\partial
x_{j}}+Au+\mu ^{2}u=0,  \tag{5.1}
\end{equation}%
\begin{equation}
\text{ }u\left( x_{1},0\right) =\psi \left( x_{1}\right) \text{, }x\in
R_{+}^{2},  \tag{5.2}
\end{equation}%
where $A=A\left( g\right) \left( x_{0},0\right) ,$ $a_{ij}=a_{ij}\left(
g\right) \left( x_{0},0\right) $ and $a_{ij}\left( g\right) $ are defined by 
$\left( 2.7\right) $ and $\psi \in h^{2,\alpha }\left( A\right) $. By
applying the Fourier transform to the problem $\left( 5.1\right) -\left(
5.2\right) $ with respect to $x_{1}$ we get%
\begin{equation}
-a_{22}\frac{d^{2}\hat{u}}{dy^{2}}+2a_{12}i\eta \frac{d\hat{u}}{dy}+\left(
A+\eta ^{2}+\mu ^{2}\right) \hat{u}=0,  \tag{5.3}
\end{equation}

\begin{equation}
\hat{u}\left( \eta ,0\right) =\hat{\psi}\left( \eta \right) ,\text{ }\eta
\in \mathbb{R},\text{ }y\in \mathbb{R}_{+},  \tag{5.4}
\end{equation}%
where $x_{2}$ and $\hat{u}\left( \eta ,x_{2}\right) $ are denoted by $y$ and 
$\hat{u}=\hat{u}\left( \eta ,y\right) ,$ respectively$.$

Let 
\begin{equation*}
p\left( \xi \right) =\xi _{1}^{2}+2a_{12}\xi _{1}\xi _{2}+a_{22}\xi _{2}^{2}%
\text{ for }\xi =\left( \xi _{1},\xi _{2}\right) \in R^{2}.
\end{equation*}%
\textbf{\ }There exists an $\alpha $ $>0$ with%
\begin{equation}
p\left( \xi \right) \geq \alpha \left\vert \xi \right\vert ^{2}\text{ for
all }\xi \in R^{2}.  \tag{5.5}
\end{equation}

The condition $\left( 5.5\right) $ implies that%
\begin{equation}
a_{12}^{2}-a_{22}\geq \alpha .  \tag{5.6}
\end{equation}

Moreover, we define%
\begin{equation}
q_{\eta }\left( \mu ,\lambda \right) =-a_{22}\lambda ^{2}+2a_{12}i\eta
\lambda +\mu ^{2}+\eta ^{2}+1=0.  \tag{5.7}
\end{equation}

\textbf{Remark 5.1. }Given $\eta \in \mathbb{R}$ and $z\in \mathbb{C},$ then
in view of $\left( 5.5\right) $ there is exactly one root of the equation $%
\left( 5.7\right) $ with positive real part. It is given by%
\begin{equation}
\lambda \left( \eta ,\mu \right) =ia\left( \eta \right) +b\left( \eta ,\mu
\right) ,\text{ }  \tag{5.8}
\end{equation}%
where 
\begin{equation*}
a\left( \eta \right) =-\frac{a_{12}}{a_{22}}\eta \text{, }b\left( \eta ,\mu
\right) =\frac{1}{a_{22}}\left[ a_{22}\left( 1+\mu ^{2}\right) +\left(
a_{12}^{2}-a_{22}\right) \eta ^{2}\right] ^{\frac{1}{2}}.
\end{equation*}%
We put 
\begin{equation*}
\tilde{X}=h^{\alpha }\left( R_{+}^{2};E\right) \text{, }\tilde{Y}%
=h^{2,\alpha }\left( R_{+}^{2};E\left( A\right) ;E\right) \text{, }
\end{equation*}%
\begin{equation*}
N_{\mu }\left( \eta ,y\right) =\exp \left\{ -\lambda \left( \eta ,\mu
\right) A_{\mu }^{\frac{1}{2}}\left( \eta \right) y\right\} .
\end{equation*}%
The main result of this section is the following:

\textbf{Theorem 5.1. }Assume the Condition 3.1 is satisfied. Then problem $%
\left( 5.1\right) -\left( 5.2\right) $ has a unique solution $u\in \tilde{Y}$
for $\psi \in h^{2,\alpha }\left( A\right) $ and $u$ is represented by%
\begin{equation*}
u\left( x_{1},x_{2}\right) =\left( \tilde{K}_{0}+\mu ^{2}\right) \psi
=F^{-1}N_{\mu }\left( \eta ,x_{2}\right) F\psi \text{ }\left( \eta \right) .
\end{equation*}

Moreover, the estimate holds%
\begin{equation}
\dsum\limits_{j=0}^{2}\left\vert \mu \right\vert ^{2-j}\left\Vert \frac{%
\partial ^{j}u}{\partial x_{1}^{j}}\right\Vert _{\tilde{X}%
}+\dsum\limits_{i,j=1}^{2}\left\Vert \frac{\partial ^{2}u}{\partial
x_{i}\partial x_{j}}\right\Vert _{\tilde{X}}+\left\Vert Au\right\Vert _{%
\tilde{X}}\leq C\left\Vert \psi \right\Vert _{h^{2,\alpha }\left( A\right) }.
\tag{5.9}
\end{equation}

For the proof we need some preparation. Here, 
\begin{equation*}
A_{\mu }=A_{\mu }\left( \eta \right) =A+\mu ^{2}+\eta ^{2}\text{, }%
q_{0}\left( \eta ,y\right) =AN_{\mu }\left( \eta ,y\right) \text{, }
\end{equation*}

\begin{equation*}
q_{j}\left( \eta ,y\right) =\left\vert \mu \right\vert ^{2-j}\eta ^{j}N_{\mu
}\left( \eta ,y\right) .
\end{equation*}%
We need the following lemmas:

\textbf{Lemma 5.1. }Assume the Condition 3.1 is satisfied. Then there exists
a unique solution $\hat{u}\left( \eta ,y\right) $ of $\left( 5.3\right)
-\left( 5.4\right) $ expressing as 
\begin{equation}
\hat{u}\left( \eta ,y\right) =N_{\mu }\left( \eta ,y\right) \hat{\psi}\left(
\eta \right)  \tag{5.10}
\end{equation}%
Moreover, the following uniformly in $\eta $ estimate holds%
\begin{equation}
\left\Vert A\hat{u}\right\Vert _{h^{\alpha }\left( \mathbb{R}_{+};E\right)
}+\dsum\limits_{j=0}^{2}\left\vert \mu \right\vert ^{2-j}\left\Vert \hat{u}%
^{\left( j\right) }\right\Vert _{h^{\alpha }\left( \mathbb{R}_{+};E\right)
}\leq C\left\Vert \hat{\psi}\right\Vert _{h^{2,\alpha }\left( A\right) }. 
\tag{5.11}
\end{equation}

\textbf{Proof. }In view of positivity of operator $A$ we know that $-A_{\mu
}^{\frac{1}{2}}\left( \eta \right) $ is analytic semigroup in $E$ (see e.g. $%
\left[ \text{25, \S\ 1.18}\right] $). The  the equation $\left( 5.3\right) $
have a solution $\hat{u}=U_{\eta ,\mu }\left( y\right) \hat{\psi}\left( \eta
\right) $ on $\left( 0,\infty \right) ,$ where 
\begin{equation*}
U_{\eta ,\mu }\left( y\right) =\exp \left\{ -\lambda \left( \eta ,\mu
\right) A_{\mu }^{\frac{1}{2}}\left( \eta \right) y\right\} \hat{\psi}\left(
\eta \right) 
\end{equation*}%
and $\lambda \left( \eta ,\mu \right) $ is a root of $\left( 5.7\right) $
with positive real part. Then from the above expresion and the properties of
analytic semigroups we get the uniform estimate%
\begin{equation*}
\left\Vert A\hat{u}\right\Vert _{h^{\alpha }\left( \mathbb{R}_{+};E\right)
}+\dsum\limits_{j=0}^{2}\left\vert \mu \right\vert ^{2-j}\left\Vert \hat{u}%
^{\left( j\right) }\right\Vert _{h^{\alpha }\left( \mathbb{R}_{+};E\right)
}\leq 
\end{equation*}

\begin{equation*}
C_{1}\sup\limits_{\eta \in \left[ 0,1\right] }\left\Vert \eta ^{1-\alpha
}AU\left( \eta \right) \hat{\psi}\left( .,\eta \right) \right\Vert
_{h^{\alpha }\left( \mathbb{R}_{+};E\right) }\leq C\left\Vert \hat{\psi}%
\right\Vert _{h^{2,\alpha }\left( A\right) },
\end{equation*}%
where $U\left( y\right) $ is a semigroup generated by $-A.$

\textbf{Lemma 5.2. }Assume the Condition 3.1 is satisfied. Then operator
functions $q_{0}\left( \eta ,y\right) $ and $q_{j}\left( \eta ,y\right) $
are Fourier multipliears in $h^{\alpha }\left( \mathbb{R};E\right) \ $%
uniformly with respect to $y\in \mathbb{R}_{+}.$

\textbf{Proof. } In view of $\left( 5.6\right) $ and the Remark 5.1 we get
that 
\begin{equation*}
\sup\limits_{\eta \in \mathbb{R},y\in \mathbb{R}_{+}}\left( 1+\left\vert
\eta \right\vert \right) ^{\frac{1}{2}}\left\Vert \text{ }\frac{d}{d\eta }%
\left[ q_{0}\left( \eta ,y\right) \right] \right\Vert _{B\left( E\right)
}\leq C_{1},\text{ }
\end{equation*}%
\begin{equation*}
\sup\limits_{\eta \in \mathbb{R},y\in \mathbb{R}_{+}}\left( 1+\left\vert
\eta \right\vert \right) ^{\frac{1}{2}}\left\Vert \text{ }\frac{d}{d\eta }%
\left[ q_{j}\left( \eta ,y\right) \right] \right\Vert _{B\left( E\right)
}\leq C_{2},
\end{equation*}%
Then by using the multiplier results in $\left[ \text{5}\right] $ we can
prove that the operator function $q_{0}\left( .,y\right) $, $q_{j}\left(
.,y\right) $ are multipliers in $h^{\alpha }\left( \mathbb{R};E\right) $
uniformly in $y\in \mathbb{R}_{+}$ and $\mu \in \mathbb{R}.$

Let 
\begin{equation*}
\Phi _{0}\left( y\right) =\left\Vert q_{0}\left( .,y\right) \right\Vert
_{M\left( h^{\alpha }\left( \mathbb{R};E\right) \right) },\text{ }\Phi
_{j}\left( y\right) =\left\Vert q_{j}\left( .,y\right) \right\Vert _{M\left(
h^{\alpha }\left( \mathbb{R};E\right) \right) }.
\end{equation*}

\bigskip \textbf{Lemma 5.3.} Assume the Condition 3.1 is satisfied. Then 
\begin{equation*}
\Phi _{0}\left( y\right) \rightarrow 0\text{, }\Phi _{j}\left( y\right)
\rightarrow 0\text{ as }y\rightarrow \infty .
\end{equation*}

\textbf{Proof. }Indeed, by properties of Fourier multiplier operators from
in $h^{\alpha }\left( \mathbb{R};E\right) $ and the theory of analytic
semigroups there exists $\omega >0$ such that we have 
\begin{equation*}
\left\Vert q_{0}\left( .,y\right) \right\Vert _{M\left( h^{\alpha }\left( 
\mathbb{R};E\right) \right) }\leq C_{1}\left( \left\vert \lambda \left( \eta
,\mu \right) \right\vert y\right) ^{-1}\exp \left\{ -\omega \left\vert
\lambda \left( \eta ,\mu \right) \right\vert y\right\} ,
\end{equation*}%
\begin{equation}
\left\Vert q_{j}\left( .,y\right) \right\Vert _{M\left( h^{\alpha }\left( 
\mathbb{R};E\right) \right) }\leq C_{2}\left\vert \mu \right\vert
^{2-j}\left( 1+\eta ^{2}\right) ^{\frac{1}{2}}\exp \left\{ -\omega
\left\vert \lambda \left( \eta ,\mu \right) \right\vert y\right\} .\text{ } 
\tag{5.12}
\end{equation}

By estimates $\left( 5.12\right) $ and Remark 5.1 we obtain the assertion.

\textbf{Proof of Theorem 5.1. } By Lemma 5.1 the problem $\left( 5.1\right)
-\left( 5.2\right) $ has a solution 
\begin{equation*}
u\left( x_{1},x_{2}\right) =\left( \tilde{K}_{0}+\mu ^{2}\right) \psi
=F^{-1}N_{\mu }\left( \eta ,x_{2}\right) \hat{\psi}\left( \eta \right) .
\end{equation*}%
By Lemmas 5.2, 5.3 the operator-functions $q_{0}\left( \eta ,y\right) $ and $%
q_{j}\left( \eta ,y\right) $ are Fourier multipliears in $h^{\alpha }\left( 
\mathbb{R};E\right) \ $uniformly with respect to $y\in \mathbb{R}_{+}.$ Then
from the estimate 5.11 we obtain trhe assertion.

Here $K_{0\mu }=K_{0\mu }\left( x_{0}\right) $ denotes the inverse of
operator $\tilde{K}_{0}\left( x_{0}\right) +\mu ^{2}$, i.e. 
\begin{equation*}
K_{0\mu }\left( x_{0}\right) =\left[ \tilde{K}_{0}\left( x_{0}\right) +\mu
^{2}\right] ^{-1}.
\end{equation*}

\textbf{Result 5.1. }From Theorem 5.1 we obtain that the opertor $%
u\rightarrow K_{0\mu }u$ is bounded from $h^{2,\alpha }\left( A\right) ^{%
\text{\ }}$into $\tilde{Y}$ and the following estimate holds%
\begin{equation*}
\dsum\limits_{j=0}^{2}\left\vert \mu \right\vert ^{2-j}\left\Vert \frac{%
\partial ^{j}K_{0\mu }}{\partial x_{1}^{j}}\right\Vert _{B\left( h^{2,\alpha
}\left( A\right) ,\tilde{X}\right) }+\dsum\limits_{i,j=1}^{2}\left\Vert 
\frac{\partial ^{2}K_{0\mu }}{\partial x_{i}\partial x_{j}}\right\Vert
_{B\left( h^{2,\alpha }\left( A\right) ,\tilde{X}\right) }+
\end{equation*}

\begin{equation*}
\left\Vert AK_{0\mu }\right\Vert _{B\left( h^{2,\alpha }\left( A\right) ,%
\tilde{X}\right) }\leq C.
\end{equation*}%
Consider now the BVP

\begin{equation*}
B\left( x_{0}\right) u=\dsum\limits_{i,j=1}^{2}-a_{ij}\frac{\partial ^{2}u}{%
\partial x_{i}\partial x_{j}}+Au+\mu ^{2}u=0,\text{ }x\in R_{+}^{2},
\end{equation*}%
\begin{equation*}
B_{0}\left( x_{0}\right) u=\left[ b_{1}\text{ }\frac{\partial }{\partial
x_{1}}u\left( x\right) +b_{2}\text{ }\frac{\partial }{\partial x_{2}}u\left(
x\right) \right] \mid _{x_{2=0}}=\psi \left( x_{1}\right) ,
\end{equation*}
where 
\begin{equation}
A=A\left( g\right) \left( x_{0},0\right) ,a_{ij}=a_{ij}\left( g\right)
\left( x_{0},0\right) ,b_{i}=b_{ij}\left( g\right) \left( x_{0},0\right) 
\tag{5.13}
\end{equation}%
and $a_{ij}\left( g\right) $ are defined by $\left( 2.7\right) $.

Here, 
\begin{equation*}
a_{1}\left( \eta ,\mu \right) =ib_{1}\left( x_{0}\right) \eta -b_{2}\left(
x_{0}\right) \lambda \left( \eta ,\mu \right) .
\end{equation*}

\begin{center}
\textbf{6. Regularity properties of abstract elliptic operator with constant
coefficients}
\end{center}

Consider the operator 
\begin{equation*}
u\rightarrow O_{1}u=u\rightarrow O_{1}\left( g\right) u=B_{0}\left( g\right)
K\left( g\right) u
\end{equation*}%
Here, $O_{10}$ denotes\ the constant coefficients version of $O_{1}$ fixed
in $\left( x_{0},0\right) $ by $\left( 5.13\right) ,$ i.e. 
\begin{equation*}
u\rightarrow O_{10}u=O_{1}\left( x_{0},0\right) =B_{0}K\left( g\right)
\left( x_{0},0\right) u.
\end{equation*}

Here, 
\begin{equation*}
O_{1\mu }u=\left( O_{1}+\mu ^{2}\right) u\text{, }O_{10\mu }u=\left(
O_{10}+\mu ^{2}\right) u\text{, }
\end{equation*}

\begin{equation*}
K_{_{0}}=K\left( g\left( x_{0},0\right) \right) \text{, \ }K_{\mu }\left(
g\right) =\left[ \tilde{K}\left( g\right) +\mu ^{2}\right] ^{-1}.
\end{equation*}%
We show the following result:

\textbf{Theorem 6.1. }Assume the Condition 3.1 is satisfied. Then the
operator $u\rightarrow \left( O_{10}+\mu ^{2}\right) u$ is an isomorphisim
from $h^{2,\alpha }\left( A\right) $ onto $h^{1,\alpha }\left( A\right) .$
Moreover, the following uniform estimates hold 
\begin{equation}
\dsum\limits_{i=0}^{1}\left\vert \mu \right\vert ^{2-i}\left\Vert \frac{%
\partial ^{i}u}{\partial x_{1}^{i}}\right\Vert _{h^{1,\alpha }\left(
A\right) }+\left\Vert Au\right\Vert _{h^{1,\alpha }\left( A\right) }\leq
C\left\Vert \left( O_{10}+\mu ^{2}\right) u\right\Vert _{h^{2,\alpha }\left(
A\right) },  \tag{6.1}
\end{equation}

\begin{equation*}
C_{1}\left\Vert u\right\Vert _{h^{2,\alpha }\left( A\right) }\leq \left\Vert
O_{20}+\mu _{0}^{2}u\right\Vert _{h^{1,\alpha }\left( A\right) }\leq
C_{2}\left\Vert u\right\Vert _{h^{2,\alpha }\left( A\right) }
\end{equation*}%
for all $u\in h^{2,\alpha }\left( A\right) $ and $\mu >0$ with sufficiently
large $\mu _{0}.$

\textbf{Proof. }\ In view of \ Lemma 3.1 and Theorem 5.1 it is clear to see
that the solution of the problem%
\begin{equation*}
\left( O_{10}+\mu ^{2}\right) u=\psi
\end{equation*}%
has a unique solution $u$ exspressed as%
\begin{equation*}
u=F^{-1}a_{1}\left( \eta ,\mu \right) N_{\mu }\left( \eta ,x_{2}\right)
F\psi \left( \eta \right) .
\end{equation*}%
By lemmas 5.1-5.3 and by multiplier result in $\left[ 5\right] $ we get that
the operator functions $\left\vert \mu \right\vert ^{2-i}\eta ^{i}K\left(
\eta ,x_{2}\right) $ and $AK\left( \eta ,x_{2}\right) $ are multipliers in $%
h^{\alpha }\left( \mathbb{R};E\right) \ $uniformly with respect to $x_{2}\in 
\mathbb{R}_{+},$ where 
\begin{equation*}
K\left( \eta ,x_{2}\right) =a_{1}\left( \eta ,\mu \right) N_{\mu }\left(
\eta ,x_{2}\right) .
\end{equation*}%
Hence, we obtain the estimate 
\begin{equation*}
\dsum\limits_{i=0}^{1}\left\vert \mu \right\vert ^{2-i}\left\Vert \frac{%
\partial ^{i}u}{\partial x_{1}^{i}}\right\Vert _{h^{1,\alpha }\left(
A\right) }+\left\Vert A_{0}u\right\Vert _{h^{1,\alpha }\left( A\right) }\leq
C\left\Vert \left( O_{10}+\mu ^{2}\right) u\right\Vert _{h^{2,\alpha }\left(
A\right) }.
\end{equation*}%
Moreover, by definitions of the\ space $h^{2,\alpha }\left( A\right) $ we
get $\left( 6.1\right) $ and the estimate%
\begin{equation*}
\left\Vert O_{10}+\mu _{0}^{2}u\right\Vert _{h^{1,\alpha }\left( A\right)
}\leq C_{2}\left\Vert u\right\Vert _{h^{2,\alpha }\left( A\right) }.
\end{equation*}%
\ Then in view of the above estimate, by reasoning as in the Theorem 5.1 we
get the assertion and corresponding esimates.

By reasoning as in Theorem 6.1 we obtain

\textbf{Theorem 6.2. }Assume the Condition 3.1 is satisfied. Then the
operator $O_{10}$ is positive and $-O_{10}$ is a generator of an analytic
semigroup in $h^{1,\alpha }\left( A\right) $.

\textbf{Proof}. \ Indeed, for positivity of the operator $O_{10}$ in $%
h^{1,\alpha }\left( A\right) $ we need to show the estimate 
\begin{equation*}
\left\Vert \left( O_{10}+\mu ^{2}\right) ^{-1}\right\Vert _{B\left(
h^{1,\alpha }\left( A\right) \right) }\leq C\mu ^{-2},\text{ }
\end{equation*}%
i.e. we have to prove the estimate%
\begin{equation*}
\left\Vert u\right\Vert _{h^{2,\alpha }\left( A\right) ^{\text{\ \ }}}\leq
C\mu \left\Vert \left( O_{10}+\mu ^{2}\right) ^{-1}u\right\Vert
_{h^{1,\alpha }\left( A\right) }
\end{equation*}%
for $u\in h^{2,\alpha }\left( A\right) $. By reasoning as the above we get
that the function $\mu ^{2}\eta a_{2}\left( \eta ,\mu \right) $ is a
multiplier in $h^{\alpha }\left( \mathbb{R};E\right) \ $uniformly with
respect to $x_{2}\in \mathbb{R}_{+},$ where 
\begin{equation*}
a_{2}\left( \eta ,x_{2}\right) =a_{1}\left( \eta ,\mu \right) N_{\mu }\left(
\eta ,x_{2}\right) .
\end{equation*}%
So, the operator $O_{10}$ is positive in $h^{1,\alpha }\left( A\right) $.
Then $-O_{10}$ is a generator of an analytic semigroup in $h^{1,\alpha
}\left( A\right) $.

Here 
\begin{equation}
Q_{2}u=O_{2}\left( g\right) u=\partial B_{0}\left( g\right) \left[ u,K\left(
g\right) g\right]  \tag{6.2}
\end{equation}%
for $g\in h_{+}^{2+\alpha }$ and $u\in h^{2,\alpha }\left( A\right) .$

Consider the operator $O_{20}$ that is a constant coefficients version of $%
O_{2}$ with fixed in $\left( x_{0},0\right) $ by $\left( 5.15\right) $, i.e. 
\begin{equation}
O_{20}=O_{20}\left( g\right) =\partial B_{0}\left( g\right) \left[ .,K\left(
g\right) g\right] \left( x_{0},0\right)  \tag{6.3}
\end{equation}%
for $g\in h^{2,\alpha }\left( A\right) .$

We prove the following result:

\textbf{Theorem 6.3. }Assume the Conditions 3.1 and 3.2 are satisfied. Then
the operator $u\rightarrow \left( O_{20}+\mu ^{2}\right) u$ is an
isomorphisim\ from $h^{2,\alpha }\left( A\right) $ onto $h^{1,\alpha }\left(
A\right) .$ Moreover, the following estimate holds%
\begin{equation}
\dsum\limits_{i=0}^{1}\left\vert \mu \right\vert ^{2-i}\left\Vert \frac{%
\partial ^{i}u}{\partial x_{1}^{i}}\right\Vert _{h^{1,\alpha }\left(
A\right) }+\left\Vert A_{0}u\right\Vert _{h^{1,\alpha }\left( A\right) }\leq
C\left\Vert \left( O_{20}+\mu ^{2}\right) u\right\Vert _{h^{2,\alpha }\left(
A\right) },  \tag{6.4}
\end{equation}

\begin{equation*}
C_{1}\left\Vert u\right\Vert _{h^{2,\alpha }\left( A\right) }\leq \left\Vert
O_{20}+\mu _{0}^{2}u\right\Vert _{h^{1,\alpha }\left( A\right) }\leq
C_{2}\left\Vert u\right\Vert _{h^{2,\alpha }\left( A\right) }
\end{equation*}%
for all $u\in h^{2,\alpha }\left( A\right) $, for sufficiently large $\mu
_{0}$ and $\mu >0.$

\textbf{Proof. }By Lemma 3.1 we have%
\begin{equation*}
O_{20}u=-\frac{\partial \upsilon }{\partial x_{1}}u_{x_{1}}+\frac{\partial
\upsilon }{\partial x_{2}}\left( \Lambda _{1}u+\Lambda _{2}u\right) ,
\end{equation*}%
where 
\begin{equation}
\upsilon =\upsilon \left( g\right) =K\left( g\right) g,\text{ }\Lambda _{1}u=%
\frac{\left( 1+g_{x_{1}}^{2}\right) u}{\nu +g},\text{ }\Lambda _{2}u=-\frac{1%
}{\nu +g}2g_{x_{1}}u_{x_{1}},  \tag{6.5}
\end{equation}%
\ 
\begin{equation*}
O_{20}u=-\frac{\partial \upsilon _{0}}{\partial x_{1}}u_{x_{1}}+\frac{%
\partial \upsilon _{0}}{\partial x_{2}}\left( \Lambda _{1}u+\Lambda
_{2}u\right) ,
\end{equation*}%
\begin{equation*}
\upsilon _{0}=\upsilon _{0}\left( g\right) =K\left( g\right) \left(
x_{0},0\right) g,\text{ }\Lambda _{1}u=\frac{\left( 1+g_{x_{1}}^{2}\right) u%
}{\nu +g},\text{ }\Lambda _{2}u=-\frac{1}{\nu +g}2g_{x_{1}}u_{x_{1}},
\end{equation*}%
By using lemmas 5.1-5.3\ we get that the operator function 
\begin{equation*}
\left( \eta ^{2}+\mu ^{2}+A_{0}\right) \left[ \eta +\Lambda \left( u\right)
a_{1}\left( \eta ,\mu \right) \right] N_{0\mu }\left( \eta ,x_{2}\right)
\end{equation*}%
is a multiplier in $h^{\alpha }\left( \mathbb{R};E\right) $ uniformly with
respect to $x_{2}\in \mathbb{R}_{+},$ where 
\begin{equation}
A_{0\mu }=A_{0\mu }\left( \eta \right) =A+\mu ^{2}+\eta ^{2},\text{ }N_{0\mu
}\left( \eta ,y\right) =A_{0\mu }^{-\frac{1}{2}}\left( \eta \right) \exp
\left\{ -\lambda \left( \eta ,\mu \right) A_{0\mu }^{\frac{1}{2}}\left( \eta
\right) y\right\} .  \tag{6.6}
\end{equation}%
Then by reasoning as in the Theorem 6.1 we obtain the assertion.

From Theorem 6.3 we obtain

\textbf{Result 6.1. }Assume the conditions 3.1 and 3.2 are satisfied. Then
the operator $O_{20}$ is positive and $-O_{20}$ is a generator of an
analityc semigroup in $h^{1,\alpha }\left( A\right) $.

Consider first of all, the BVP 
\begin{equation}
-\dsum\limits_{k,j=1}^{2}a_{kj}\frac{\partial ^{2}u}{\partial x_{k}\partial
x_{j}}+\left( A+\mu ^{2}\right) u+\mu ^{2}u=V\left( x\right) ,  \tag{6.7}
\end{equation}%
\begin{equation}
\text{ }u\left( x_{1},0\right) =0\text{, }x=\left( x_{1},x_{2}\right) \in
R_{+}^{2},  \tag{6.8}
\end{equation}%
where $A=A\left( g\right) \left( x_{0},0\right) ,$ $a_{kj}=a_{kj}\left(
g\right) \left( x_{0},0\right) $ and $a_{kj}\left( g\right) $ are defined by 
$\left( 5.15\right) $. \textbf{\ }

Let $S_{0}=S\left( g\right) \left( x_{0},0\right) $ denotes the realizasion
operator in $\tilde{X}$ generated by $\left( 6.7\right) -\left( 6.8\right) $
for $\mu =0$, i.e.%
\begin{equation*}
D\left( S_{0}\right) =\tilde{Y}\text{, }S_{0}u=-\dsum%
\limits_{k,j=1}^{2}a_{kj}\frac{\partial ^{2}u}{\partial x_{k}\partial x_{j}}%
+Au.
\end{equation*}

From $\left[ \text{4, Theorem 2}\right] $ we obtain the following:

\bigskip \textbf{Result 6.2. }Assume the Condition 3.1 is satisfied. Then;

(1) problem $\left( 6.7\right) -\left( 6.8\right) $ for sufficientli large $%
\mu >0$ has a unique solution $u\in \tilde{Y}$ for $V\in \tilde{X};$

(2) the uniform coercive estimate holds%
\begin{equation}
\dsum\limits_{i=0}^{1}\left\vert \mu \right\vert ^{2-i}\left\Vert \frac{%
\partial ^{i}u}{\partial x_{1}^{i}}\right\Vert _{\tilde{X}%
}+\dsum\limits_{k,j=1}^{2}\left\Vert \frac{\partial ^{2}u}{\partial
x_{k}\partial x_{j}}\right\Vert _{\tilde{X}}+\left\Vert Au\right\Vert
_{X}\leq C\left\Vert V\right\Vert _{\tilde{X}};  \tag{6.9}
\end{equation}

(3) the operator $S_{0}$ is a positive and $-S_{0}$ is a generator of an
analytic semigrop in $h^{1,\alpha }\left( A\right) .$

The estimate $\left( 6.9\right) $ particularly, implies that $\left(
S_{0}+\mu ^{2}\right) ^{-1}\in B\left( \tilde{X},\tilde{Y}\right) .$

Consider the inhomogenous problem 
\begin{equation}
-\dsum\limits_{k,j=1}^{2}a_{ij}\frac{\partial ^{2}u}{\partial x_{k}\partial
x_{j}}+Au+\mu ^{2}u=V\left( x\right) ,  \tag{6.10}
\end{equation}%
\begin{equation}
\gamma u=\text{ }u\left( x_{1},0\right) =\psi \left( x_{1}\right) \text{, }%
x=\left( x_{1},x_{2}\right) \in R_{+}^{2}.  \tag{6.11}
\end{equation}

\textbf{Theorem 6.4.} Assume the conditions 3.1 and 3.2 are satisfied. Then
the operator $u\rightarrow G_{0}u=\left\{ L_{0}u,\text{ }\gamma u\right\} $
is an isomorphisim from $\tilde{Y}$ onto $\tilde{X}\times h^{2,\alpha
}\left( A\right) .$

\textbf{Proof. }From definition of $\tilde{X}$, $\tilde{Y}$, $h^{2,\alpha
}\left( A\right) $, from expresion of $L_{0}$ and by virtue of trace result
in $\tilde{Y}$ $\left[ \text{15, Ch.2}\right] $ we get that 
\begin{equation*}
\left\Vert G_{0}u\right\Vert _{\tilde{X}\times h^{1,\alpha }\left( A\right)
}=\left\Vert L_{0}u\right\Vert _{\tilde{X}}+\left\Vert \gamma u\right\Vert
_{h^{2,\alpha }\left( A\right) }\leq C\left\Vert u\right\Vert _{\tilde{Y}},
\end{equation*}%
i.e. the operator $G_{0}$ is bounded linear from $\tilde{Y}$ into $\tilde{X}%
\times h^{2,\alpha }\left( A\right) .$ Hence, in view of Banach theorem it
is sufficient to show that the operator $G_{0}$ is inective and surjective
from $\tilde{Y}$ onto $\tilde{X}\times h^{2,\alpha }\left( A\right) $. From
Theorem 5.1 we obtain that the corresponding homogenous problem $L_{0}u=0,$ $%
\gamma u=0$ has a zero solution, i.e.\ the operator $G_{0}$ is inective. So,
it remain to show that this operator is surjective. By Theorem 5.1 we obtain
that problem $L_{0}u=0,$ $\gamma u=\psi \ $has a solution $u_{1}\in $ $%
\tilde{Y}$ for all $\psi \in h^{2,\alpha }\left( A\right) .$ Moreover, from
the Result 6.2 we get that problem $L_{0}u=V,$ $\gamma u=0$ has a solution $%
u_{2}\in $ $\tilde{Y}$ for all $V\in \tilde{X}.$ Then $u=u_{1}+u_{2}$ is a
solution of $\left( 6.10\right) -\left( 6.11\right) $ that belongs to $%
\tilde{Y},$ i.e. the operator $G_{0}$ is surgective from $\tilde{Y}$ onto $%
\tilde{X}\times h^{2,\alpha }\left( A\right) .$

From Theorems 5.1 and 6.4 we obtain the following

\textbf{Result 6.3. }The soluion $u$ of the problem $\left( 6.10\right)
-\left( 6.11\right) $ is exspressed as 
\begin{equation*}
u\left( x\right) =S_{1,\mu }V+S_{2,\mu }\left( \psi -\gamma u_{1}\right) ,
\end{equation*}%
where 
\begin{equation*}
S_{1,\mu }V=r_{+}F^{-1}\left( p\left( \xi \right) +A+\mu ^{2}\right) ^{-1}F%
\tilde{V},\text{ }S_{2,\mu }\upsilon =F^{-1}N_{\mu }\left( \xi
_{1},x_{2}\right) F\upsilon ,
\end{equation*}%
here $r_{+}$ is the restriction operator from $R^{2}$ into $R_{+}^{2}$ and $%
\tilde{V}=\tilde{V}\left( x_{1},x_{2}\right) $ is an exstension of $V\left(
x_{1},x_{2}\right) $ on $R^{2}$, i.e 
\begin{equation*}
\tilde{V}\left( x_{1},x_{2}\right) =\left\{ 
\begin{array}{c}
V\left( x_{1},x_{2}\right) \text{, if }x_{1},x_{2}\in \bar{R}_{+}^{2} \\ 
V\left( x_{1},-x_{2}\right) \text{ if }x_{1},x_{2}\in R_{+}^{2}%
\end{array}%
\right.
\end{equation*}

$N_{\mu }\left( \xi _{1},x_{2}\right) $-is operator function defined by $%
\left( 5.8\right) $\ and $u_{1}$ is a solution of the equation 
\begin{equation*}
-\dsum\limits_{k,j=1}^{2}a_{ij}\frac{\partial ^{2}u}{\partial x_{k}\partial
x_{j}}+Au+\mu ^{2}u=\tilde{V}\left( x\right) ,\text{ }x\in R^{2}.
\end{equation*}

Let $O_{3}=O_{3}\left( g\right) $ be the operator in $\left( 4.6\right) $
defined by 
\begin{equation}
O_{3}u=B_{0}\left( g\right) S\left( g\right) \partial B\left( g\right) \left[
u,K\left( g\right) g\right]  \tag{6.12}
\end{equation}%
for $g\in h^{2,\alpha }\left( A\right) $ and $u\in Y.$

In view of Lemma 3.1 and Lemma 3.2 we get 
\begin{equation*}
O_{3}u=B_{0}\left( g\right) S\left( g\right) \partial B\left( g\right) \left[
u,\upsilon \right] ,\text{ }\upsilon =\upsilon _{g}\left( x\right) =\left(
K\left( g\right) g\right) \left( x\right) ,
\end{equation*}%
where%
\begin{equation*}
\partial B\left( g\right) \left[ u,\upsilon \right] =\frac{2\beta }{\nu +g}%
\left\{ \left[ \left( \frac{g_{x_{1}}u}{\nu +g}-u_{x_{1}}\right) \right]
\upsilon _{x_{1}x_{2}}+\right.
\end{equation*}%
\begin{equation}
\frac{1}{\left( \nu +g\right) ^{2}}\left[ \left( \frac{1}{\beta }+\beta
g_{x_{1}}^{2}\right) u-\frac{\beta }{\nu +g}g_{x_{1}}u_{x_{1}}\right]
\upsilon _{x_{2}x_{2}}+\partial A\left( g\right) \left[ u,\upsilon \right] -%
\text{ }  \tag{6.13}
\end{equation}

\begin{equation*}
\left. \left[ \frac{g_{x_{1}}^{2}u}{\left( \nu +g\right) ^{2}}-\frac{%
g_{x_{1}x_{1}}u+4g_{x_{1}}u_{x_{1}}}{2\left( \nu +g\right) }+\frac{%
u_{x_{1}x_{1}}}{2}\right] \upsilon _{x_{2}}\right\} ,
\end{equation*}

\begin{equation*}
O_{3}=O_{3}\left( g\right) =\dsum\limits_{i=1}^{4}O_{3i}\left( g\right) u,%
\text{ }O_{3i}\left( g\right) u=B_{0}\left( g\right) S\left( g\right)
G_{i}\left( g\right) ,
\end{equation*}

where 
\begin{equation*}
G_{i}\left( g\right) =G_{i}\left( g\right) \left[ u,\upsilon \right] ,\text{ 
}G_{1}\left( g\right) =\frac{2\beta }{\nu +g}\left[ \left( \frac{g_{x_{1}}u}{%
\nu +g}-u_{x_{1}}\right) \right] \upsilon _{x_{1}x_{2}},
\end{equation*}

\begin{equation}
G_{2}\left( g\right) =\left[ \frac{u}{\left( \nu +g\right) ^{2}}\left( \frac{%
1}{\beta }+\beta g_{x}^{2}\right) -\frac{\beta }{\nu +g}g_{x_{1}}u_{x_{1}}%
\right] \upsilon _{x_{2}x_{2}},\text{ }  \tag{6.14}
\end{equation}%
\begin{equation*}
G_{3}\left( g\right) =\frac{2\beta }{\nu +g}\partial A\left( g\right) \left[
u,\upsilon \right] \text{, }
\end{equation*}%
\begin{equation*}
G_{4}\left( g\right) =-\left[ \left( \frac{g_{x_{1}}^{2}}{\left( \nu
+g\right) ^{2}}-\frac{g_{x_{1}x_{1}}}{2\left( \nu +g\right) }\right) u-\frac{%
4g_{x_{1}}u_{x_{1}}}{2\left( \nu +g\right) }+\frac{u_{x_{1}x_{1}}}{2}\right]
\upsilon _{x_{2}}.
\end{equation*}%
Consider the operator $O_{30}$ that is the constant coefficients version of $%
O_{3}$ fixed in $\left( x_{0},0\right) $ which defined by $\left(
5.16\right) $, i.e.\ from the above equality and from $\left( 6.13\right) $
we get $O_{30k}=B_{0}\left( g\right) S\left( g\right) G_{k}\left( g\right)
\left( x_{0},0\right) ,$%
\begin{equation}
O_{30}u=O_{30}\left( g\right) u=O_{301}u+O_{302}u+O_{303}u+O_{304}u, 
\tag{6.15}
\end{equation}%
where%
\begin{equation*}
O_{30i}=O_{3i}\left( g\left( x_{0},0\right) \right) u
\end{equation*}

\begin{equation*}
w\left( x\right) =w_{g}\left( x\right) =\frac{\beta }{\nu +g}\frac{\partial 
}{\partial x_{2}}\upsilon _{g}\left( x_{1},x_{2}\right) \text{, }\left(
x_{1},x_{2}\right) \in \Omega ,\text{ }w_{0}=w_{g}\left( x^{0},0\right)
\end{equation*}%
For $u\in h^{2,\alpha }\left( A\right) $ define the operator by%
\begin{equation*}
P_{1}u\left( x_{1},x_{2}\right) =\left( P_{1}\left( g,x_{1}^{0}\right)
u\right) \left( x_{1},x_{2}\right) =w_{g}\left( x\right)
u_{x_{1}x_{1}}e^{-x_{2}}.
\end{equation*}

For later purposes we need the following technical lemmas:

\textbf{Lemma 6.1.} Assume the conditions 3.1 and 3.2 are satisfied. Then:

(a) $P_{1}\in B\left( h^{2,\alpha }\left( A\right) ,\tilde{X}\right) ;$

(b) There exists a positive constant $C=C(g)$ such that 
\begin{equation*}
\left\Vert P_{1}u-wu_{x_{1}x_{1}}\right\Vert _{h^{\alpha }\left( \tilde{U}%
_{r};E\right) }\leq Cr\left\Vert u\right\Vert _{h^{2,\alpha }\left( A\right)
}
\end{equation*}

for all $u\in h^{2,\alpha }\left( A\right) $, $r\in \left( 0,1\right) ,$
where 
\begin{equation*}
\tilde{U}_{r}=\mathbb{R}_{+}^{2}\cap U_{r}\left( x_{1}^{0},0\right) ,
\end{equation*}

here, $U_{r}\left( x_{1}^{0},0\right) $ denotes the two dimensional ball
with radius $r$ centered at $\left( x_{1}^{0},0\right) .$

\textbf{Proof. \ }Indeed, from the expression $\left( 6.14\right) $ in view
of the Theorem 5.1 and by virtue of trace theorem in $Y$ we get (a); Then by
using the integral mean value theorem and the trace theorem we obtain%
\begin{equation*}
\left\Vert P_{1}u-wu_{x_{1}x_{1}}\right\Vert _{h^{\alpha }\left( \tilde{U}%
_{r};E\right) }\leq \left\Vert \left( we^{-x_{2}}-w\right)
u_{x_{1}x_{1}}\right\Vert _{h^{\alpha }\left( \tilde{U}_{r};E\right) }\leq
Cr\left\Vert u\right\Vert _{h^{2,\alpha }\left( A\right) }.
\end{equation*}

Let 
\begin{equation*}
\partial A=\partial A\left( x\right) =\partial A\left( g\right) \left(
x\right) ,\text{ }A_{0}=\partial A\left( g\right) \left( x_{10},0\right) .
\end{equation*}

\textbf{Lemma 6.2.} Assume the conditions 3.1 and 3.2 are satisfied. Then
the operator $u\rightarrow \left( O_{304}+\mu ^{2}\right) u$ is an
isomorphisim\ from $h^{2,\alpha }\left( A\right) $ onto $h^{1,\alpha }\left(
A\right) .$ Moreover, the following estimate holds%
\begin{equation*}
\dsum\limits_{i=0}^{1}\left\vert \mu \right\vert ^{2-i}\left\Vert \frac{%
\partial ^{i}u}{\partial x_{1}^{i}}\right\Vert _{h^{1,\alpha }\left(
A\right) }+\left\Vert A_{0}u\right\Vert _{h^{1,\alpha }\left( A\right) }\leq
C\left\Vert \left( O_{304}+\mu ^{2}\right) u\right\Vert _{h^{2,\alpha
}\left( A\right) },
\end{equation*}

\begin{equation}
C_{1}\left\Vert u\right\Vert _{h^{2,\alpha }\left( A\right) }\leq \left\Vert
O_{304}+\mu _{0}^{2}u\right\Vert _{h^{1,\alpha }\left( A\right) }\leq
C_{2}\left\Vert u\right\Vert _{h^{2,\alpha }\left( A\right) }  \tag{6.16}
\end{equation}%
for all $u\in h^{2,\alpha }\left( A\right) $ and\ for sufficiently large $%
\mu _{0}$, $\mu >0,.$

\textbf{Proof. }From the expression $\left( 6.14\right) $ and Lemma 6.1 by
reasoning as in $\left[ \text{9, lemma 5.3}\right] $ we get that the
operator $u\rightarrow \left( P_{1}+\mu ^{2}\right) u$ is an isomorphism
from $\tilde{Y}$ onto $\tilde{X}.$ The Theorem 6.4 implies that the operator 
$u\rightarrow \left[ S_{x_{0}}\left( g\right) +\mu ^{2}\right] G_{40}u$ is
an isomorphism from $\tilde{Y}$ onto $\tilde{Y}$ for sufficiently large $\mu
>0.$ Then in view of trace theorem in $\tilde{Y}$ we get that the operator $%
u\rightarrow \left[ O_{304}+\mu ^{2}\right] u$ is an isomorphism from $%
\tilde{Y}$ onto $h^{1,\alpha }\left( A\right) .$ Moreover, in view of Result
6.3 by reasoning as in the Theorem 6.1 we obtain the estimates $\left(
6.16\right) $ for all $u\in h^{2,\alpha }\left( A\right) $ and\ for $\mu >0$
with sufficiently large $\mu _{0}.$

\textbf{Lemma 6.3.} Assume the conditions 3.1 and 3.2 are satisfied. Then
the operator $u\rightarrow \left( O_{303}+\mu ^{2}\right) u$ is an
isomorphisim\ from $h^{2,\alpha }\left( A\right) $ onto $h^{1,\alpha }\left(
A\right) .$ Moreover, the following estimate holds%
\begin{equation*}
\dsum\limits_{i=0}^{1}\left\vert \mu \right\vert ^{2-i}\left\Vert \frac{%
\partial ^{i}u}{\partial x_{1}^{i}}\right\Vert _{h^{1,\alpha }\left(
A\right) }+\left\Vert A_{0}u\right\Vert _{h^{1,\alpha }\left( A\right) }\leq
C\left\Vert \left( O_{303}+\mu ^{2}\right) u\right\Vert _{h^{2,\alpha
}\left( A\right) },
\end{equation*}

\begin{equation}
C_{1}\left\Vert u\right\Vert _{h^{2,\alpha }\left( A\right) }\leq \left\Vert
O_{303}+\mu _{0}^{2}u\right\Vert _{h^{1,\alpha }\left( A\right) }\leq
C_{2}\left\Vert u\right\Vert _{h^{2,\alpha }\left( A\right) }  \tag{6.17}
\end{equation}%
for all $u\in h^{2,\alpha }\left( A\right) ^{\text{\ \ }}$and\ for
sufficiently large $\mu _{0}$, $\mu >0.$

\textbf{Proof. }From the expression $\left( 6.14\right) $ by properties of
positive operators we get that the map $u\rightarrow \left( G_{30}+\mu
^{2}\right) u$ is an isomorphism from $\tilde{Y}$ onto $\tilde{X}.$ The
Theorem 6.4 implies that the operator $u\rightarrow \left[ S_{x_{0}}\left(
g\right) +\mu ^{2}\right] G_{30}u$ is an isomorphism from $\tilde{Y}$ onto $%
\tilde{Y}$ for sufficiently large $\mu >0.$ Then in view of trace theorem in 
$\tilde{Y}$ we get that the operator $u\rightarrow \left[ O_{303}+\mu ^{2}%
\right] u$ is an isomorphism from $\tilde{Y}$ onto $h^{1,\alpha }\left(
A\right) .$ Moreover, in view of Result 6.3 by reasoning as in the Theorem
6.1 we obtain the estimates $\left( 6.17\right) $ for all $u\in h^{2,\alpha
}\left( A\right) $ and\ for $\mu >0$ with sufficiently large $\mu _{0}.$

In a similar way as Lemma 6.1 and by reasoning as in $\left[ \text{9},\text{%
lemma 5.4}\right] $ we obtain

\textbf{Lemma 6.4.} Assume the conditions 3.1 and 3.2 are satisfied. Then
the operator $u\rightarrow \left( O_{30k}+\mu ^{2}\right) u$ is an
isomorphisim\ from $h^{2,\alpha }\left( A\right) $ onto $h^{1,\alpha }\left(
A\right) .$ Moreover, the following estimate holds%
\begin{equation*}
\dsum\limits_{i=0}^{1}\left\vert \mu \right\vert ^{2-i}\left\Vert \frac{%
\partial ^{i}u}{\partial x_{1}^{i}}\right\Vert _{h^{1,\alpha }\left(
A\right) }+\left\Vert A_{0}u\right\Vert _{h^{1,\alpha }\left( A\right) }\leq
C\left\Vert \left( O_{30k}+\mu ^{2}\right) u\right\Vert _{h^{2,\alpha
}\left( A\right) },
\end{equation*}

\begin{equation}
C_{1}\left\Vert u\right\Vert _{h^{2,\alpha }\left( A\right) }\leq \left\Vert
O_{30k}+\mu _{0}^{2}u\right\Vert _{h^{1,\alpha }\left( A\right) }\leq
C_{2}\left\Vert u\right\Vert _{h^{2,\alpha }\left( A\right) }\text{, }k=1,%
\text{ }2,\text{ }  \tag{6.18}
\end{equation}%
for all $u\in h^{2,\alpha }\left( A\right) $ and $\mu >0.$

\textbf{Theorem 6.5. }Assume the conditions 3.1 and 3.2 are satisfied. Then
the operator $u\rightarrow \left( Q_{30}+\mu ^{2}\right) u$ is an
isomorphisim\ from $h^{2,\alpha }\left( A\right) $ onto $h^{1,\alpha }\left(
A\right) .$ Moreover, the following estimate holds%
\begin{equation*}
\dsum\limits_{i=0}^{1}\left\vert \mu \right\vert ^{2-i}\left\Vert \frac{%
\partial ^{i}u}{\partial x_{1}^{i}}\right\Vert _{h^{1,\alpha }\left(
A\right) }+\left\Vert A_{0}u\right\Vert _{h^{1,\alpha }\left( A\right) }\leq
C\left\Vert \left( O_{30}+\mu ^{2}\right) u\right\Vert _{h^{2,\alpha }\left(
A\right) },
\end{equation*}

\begin{equation}
C_{1}\left\Vert u\right\Vert _{h^{2,\alpha }\left( A\right) }\leq \left\Vert
O_{30}+\mu _{0}^{2}u\right\Vert _{h^{1,\alpha }\left( A\right) }\leq
C_{2}\left\Vert u\right\Vert _{h^{2,\alpha }\left( A\right) }  \tag{6.19}
\end{equation}%
for all $u\in h^{2,\alpha }\left( A\right) $ and\ $\mu >0.$

\textbf{Proof. }From the expressions $\left( 6.13\right) -\left( 6.14\right) 
$ and from lemmas 6.2-6.4\ we get that the operator $u\rightarrow
K_{1x_{0}}\left( g\right) u$ is an isomorphism from $\tilde{Y}$ onto $\tilde{%
X}.$ The Theorem 6.4 implies that the operator $u\rightarrow \left[
S_{x_{0}}\left( g\right) +\mu ^{2}\right] K_{1x_{0}}\left( g\right) u$ is an
isomorphism from $\tilde{Y}$ onto $\tilde{Y}$ for sufficiently large $\mu
>0. $ Then in view of trace theorem in $\tilde{Y}$ we get that the operator $%
u\rightarrow \left[ O_{30}+\mu ^{2}\right] u$ is an isomorphism from $\tilde{%
Y}$ onto $h^{1,\alpha }\left( A\right) $ for sufficiently large $\mu _{0}$
and $\mu >0.$

From Theorem 6.5 we obtain

\textbf{Result 6.4. }Assume the Conditions 3.1 and 3.2 are satisfied. Then
the operator $O_{30}$ is positive and $-O_{30}$ is a generator of an
analytic semigroup in $h^{1,\alpha }\left( A\right) $.

From Theorems 6.1, 6.3 and 6.5 we obtain

\textbf{Result 6.5. \ }Assume the Conditions 3.1 and 3.2 are satisfied. Then
\ 
\begin{equation*}
O_{k0}+\mu _{0}\in H\left( h^{2,\alpha }\left( A\right) ,h^{1,\alpha }\left(
A\right) \right) \text{, }k=1,2,3
\end{equation*}

for sufficiently large $\mu _{0}>0.$

In view of theorems 6.3, 6.5 \ and by lemmas 2.3, 3.1 and Theorem 3.1 we
obtain

\textbf{Result 6.6. }Assume the Conditions 3.1 and 3.2 are satisfied. The
the following estimate holds%
\begin{equation*}
\left\Vert O_{20}\left( g\right) +O_{30}\left( g\right) +\mu _{0}\right\Vert
_{B\left( h^{2,\alpha }\left( A\right) ,h^{1,\alpha }\left( A\right) \right)
}\leq C\left\Vert g\right\Vert _{h^{2,\alpha }\left( A\right) ^{\text{\ \ }}}
\end{equation*}%
for all 
\begin{equation*}
g\in B_{0}\subset h^{2,\alpha }\left( A\right) ^{\text{\ \ }}.
\end{equation*}

Let 
\begin{equation*}
O_{0}=O_{01}+O_{02}+O_{03}.
\end{equation*}

Now, we will show that%
\begin{equation*}
O_{0}\in H\left( h^{2,\alpha }\left( A\right) ,h^{1,\alpha }\left( A\right)
\right) .
\end{equation*}

\textbf{Theorem 6.6. }Assume the Conditions 3.1 and 3.2 are satisfied.
Suppose 
\begin{equation*}
\left\Vert w_{0}\right\Vert _{D_{A}\left( \alpha \right) }\leq \frac{\alpha
\left( g\right) \left( x\right) }{a_{22}\left( g\right) \left( x,0\right) }=%
\frac{\alpha _{0}}{a_{22}}.
\end{equation*}%
Then the operator $u\rightarrow \left( Q_{0}+\mu ^{2}\right) u$ is an
isomorphisim\ from $h^{2,\alpha }\left( A\right) $ onto $h^{1,\alpha }\left(
A\right) .$ Moreover, the following estimate holds%
\begin{equation*}
\dsum\limits_{i=0}^{1}\left\vert \mu \right\vert ^{2-i}\left\Vert \frac{%
\partial ^{i}u}{\partial x_{1}^{i}}\right\Vert _{h^{1,\alpha }\left(
A\right) }+\left\Vert A_{0}u\right\Vert _{h^{1,\alpha }\left( A\right) }\leq
C\left\Vert \left( O_{0}+\mu ^{2}\right) u\right\Vert _{h^{2,\alpha }\left(
A\right) },
\end{equation*}

\begin{equation}
C_{1}\left\Vert u\right\Vert _{h^{2,\alpha }\left( A\right) }\leq \left\Vert
O_{0}+\mu _{0}^{2}u\right\Vert _{h^{1,\alpha }\left( A\right) }\leq
C_{2}\left\Vert u\right\Vert _{h^{2,\alpha }\left( A\right) }  \tag{6.20}
\end{equation}%
for all $u\in h^{2,\alpha }\left( A\right) $ and\ for sufficiently large $%
\mu _{0}$, $\mu >0$. Particularly, the operator $O_{0}$ is positive and $%
-O_{0}$ is a generator of an analytic semigroup in $h^{1,\alpha }\left(
A\right) $.

\textbf{Proof. }Indeed, by using Theorems 6.1, 6.3, 6.5, the Results 6.3,
6.5, by reasoning as in lemma 5.4, Theorem 5.6, Corollary 5.7 in $\left[ 2%
\right] $ and the perturbation results for space $H\left( h^{2,\alpha
}\left( A\right) ,h^{1,\alpha }\left( A\right) \right) $ we obtain the
assertion.

Let 
\begin{equation*}
W_{t}=\left\{ g\in h_{+}^{2,\alpha }\text{, }\inf\limits_{x\in \left(
-\infty ,\infty \right) }\left[ \frac{t}{\nu +g}\frac{\partial }{\partial
x_{2}}\upsilon _{g}\left( x,0\right) +k_{g}\left( x\right) >0\right]
\right\} ,
\end{equation*}

\begin{equation*}
k_{g}\left( x\right) =\frac{\alpha \left( g\right) \left( x\right) }{%
a_{22}\left( g\right) \left( x,0\right) },\text{ }t\in \left[ 0,1\right] .
\end{equation*}

\textbf{Remark 6.1. }Suppose that $g\in W_{1}.$ Then 
\begin{equation*}
w_{\pi }=w_{g}\left( x,0\right) <\frac{\alpha _{0}}{a_{22}}\text{ for }x\in
\left( -\infty ,\infty \right) .
\end{equation*}

From Theorem 5.1 we know that there is a constant $M>0$ such that 
\begin{equation*}
\left\Vert K\left( g\right) \right\Vert _{B\left( h^{2,\beta }\left( \Gamma
;H\right) ,h^{2,\beta }\left( \Omega ;H\right) \right) }\leq M\text{, }\beta
\in \left( 0,\alpha \right)
\end{equation*}

for all $g\in h_{+}^{2,\beta }$ satisfying $\left\Vert g\right\Vert
_{h^{2,\beta }}\leq \chi .$ Now define%
\begin{equation*}
k=\frac{\left( \eta -\chi \right) ^{3}\chi }{2M\left[ 1+\left( \eta +\chi
\right) ^{2}+\chi ^{2}\right] \left( 1+\chi \right) ^{2}}.
\end{equation*}

\begin{center}
\textbf{7. Coercive estimates for the linearization}
\end{center}

It is clear that $O_{0}=O_{10}+$ $O_{20}+$ $O_{30}$ may be viewed as a
principal part of the linearization of $\partial O\left( g\right) $\ with
coefficients fixed in $\left( x^{0},0\right) $. \ Our main goal in this
section is to prove that the operator $O$ belongs to the class 
\begin{equation*}
H\left( h^{2,\alpha }\left( A\right) ,h^{1,\alpha }\left( A\right) \right) .
\end{equation*}%
We use the estimates of local operators in the preceding section to derive
coercive estimates for the linearization operator $\partial O\left( g\right) 
$. Let 
\begin{equation*}
\partial O\left( g\right) =O_{1}\left( g\right) +O_{2}\left( g\right)
+O_{3}\left( g\right) ,
\end{equation*}%
where 
\begin{equation}
O_{1}\left( g\right) =B_{0}\left( g\right) K\left( g\right) ,\text{ }%
O_{2}\left( g\right) =\partial B_{0}\left( g\right) \left[ .,K\left(
g\right) g\right] ,  \tag{7.1}
\end{equation}

\begin{equation*}
\text{ }O_{3}\left( g\right) =-B_{0}\left( g\right) S\left( g\right)
\partial B\left( g\right) \left[ .,K\left( g\right) g\right] .
\end{equation*}%
Given $g\in h_{+}^{2,\alpha }$ and $t\in \left[ 0,1\right] ,$ set 
\begin{equation*}
\partial O^{t}\left( g\right) =O_{1}\left( g\right) +t\partial B_{0}\left(
g\right) \left[ .,K\left( g\right) g\right] -tB_{0}\left( g\right) S\left(
g\right) \partial B\left( g\right) \left[ .,K\left( g\right) g\right]
\end{equation*}%
and observe that%
\begin{equation*}
\partial O^{1}\left( g\right) =\partial O\left( g\right) .
\end{equation*}

Let $\delta >0$ be given and let $\left\{ V_{j},\varphi _{j}\text{, }j\in 
\mathbb{N}\right\} $ denote $\delta $-localizasion sequence for $S=\mathbb{%
R\times }\left( -\frac{\delta }{2},\frac{\delta }{2}\right) $, the covering $%
\left\{ V_{j}\text{, }j\in \mathbb{N}\right\} $ has finite multiplicity,
diam $U_{j}$ $<\delta $, and $\left\{ V_{j},\varphi _{j}\text{, }j\in 
\mathbb{N}\right\} $ is a partition of unity on $S$ with $\dsum\limits_{j\in 
\mathbb{N}}\varphi _{j}\left( x\right) \equiv 1$. Moreover, we fix $%
x_{1j}\in \mathbb{R}$ such that $\left( x_{1j},0\right) \in V_{j},$ $j\in 
\mathbb{N}.$

Here, we will prove the following result

\textbf{Theorem 7.1. }Assume the conditions 3.1 and 3.2 are satisfied.
Suppose that $W_{0}\subset W_{1}$ is compact$,$ $\beta \in \left( 0,\alpha
\right) $ and that $k>0.$ Then there exist $\delta \in $ $(0,1]$, a $\delta $%
-localization sequence $\left\{ V_{j},\varphi _{j}\text{, }j\in \mathbb{N}%
\right\} $, $\beta \in \left( 0,\alpha \right) $ and a positive constant $%
C=C(W_{0},M,\delta )$ such that 
\begin{equation*}
\left\Vert \left[ \varphi _{j}\partial O_{t}\left( g\right) -O_{\pi }\left(
g,x_{1j}\right) \varphi _{j}\right] \upsilon \right\Vert _{h^{1,\alpha
}\left( A\right) }\leq k\left\Vert \varphi _{j}\upsilon \right\Vert
_{h^{2,\alpha }\left( A\right) }+C\left\Vert \upsilon \right\Vert
_{h^{2,\beta }\left( A\right) }
\end{equation*}%
for all $\upsilon \in h^{2,\alpha }\left( A\right) ,$ $j\in \mathbb{N},$ $%
t\in \left[ 0,1\right] $ and $g\in W_{0}.$

For proving Theorem 7.1 we need some preparation. Let 
\begin{equation*}
\text{ }A_{j}=\partial A\left( g\right) \left( x_{1j},0\right) .
\end{equation*}

\textbf{\ }Consider the following equation 
\begin{equation}
\left( O\left( g\right) +\mu \right) u=f  \tag{7.2}
\end{equation}%
for%
\begin{equation*}
f=\left[ O\left( g\right) +\mu \right] u\in B_{1p}.
\end{equation*}%
From $\left( 7.2\right) $ for $u_{j}=u\varphi _{j}$, $u\in D_{A}\left(
\alpha \right) $ we get%
\begin{equation}
\left( O\left( g\right) +\mu \right) u_{j}=\dsum\limits_{k=1}^{3}O_{k}\left(
g\right) u_{j}+\mu u_{j}=\dsum\limits_{k=1}^{3}F_{kj}+f.\varphi _{j}, 
\tag{7.3}
\end{equation}%
where 
\begin{equation*}
\text{ }F_{1j}=\gamma \left[ \left( b_{1}\frac{\partial \varphi _{j}}{%
\partial x_{1}}+b_{2}\frac{\partial \varphi _{j}}{\partial x_{2}}\right)
U_{1}\left( g\right) \right] \text{, }U_{1}\left( g\right) =K\left( g\right)
u,
\end{equation*}%
\begin{equation}
F_{2j}=\gamma \left[ -\frac{\partial \upsilon }{\partial x_{1}}\frac{%
\partial \varphi _{j}}{\partial x_{1}}+\frac{2g_{x_{1}}}{\nu +g}\frac{%
\partial \upsilon }{\partial x_{2}}\frac{\partial \varphi _{j}}{\partial
x_{1}}\right] u,\text{ }\upsilon =\upsilon \left( g\right) =K\left( g\right)
g,  \tag{7.4}
\end{equation}%
\begin{equation*}
F_{3j}=\gamma \left[ \left( b_{1}\frac{\partial \varphi _{j}}{\partial x_{1}}%
+b_{2}\frac{\partial \varphi _{j}}{\partial x_{2}}\right) U_{3}\left(
g\right) \right] ,\text{ }U_{3}\left( g\right) =S\left( g\right) \partial
B\left( g\right) \left[ u,K\left( g\right) g\right] .
\end{equation*}%
By freezing in $\left( 7.3\right) $ coefficients at points $\left(
x_{1j},0\right) $ we have localized equations 
\begin{equation}
\left( O\left( x_{1j}\right) +\mu \right) u_{j}=F_{j},  \tag{7.5}
\end{equation}%
where 
\begin{equation*}
F_{j}=\dsum\limits_{k=1}^{3}F_{kj}+f.\varphi
_{j}+\dsum\limits_{k=1}^{3}\left( O_{kj}-O_{j}\right) u_{j},
\end{equation*}%
\begin{equation*}
\left( O_{1j}-O_{1}\right) u_{j}=B_{0}\left( g\right) \left[ K_{j}-K\left(
g\right) \right] u_{j}+\left[ B_{0j}-B_{0}\left( g\right) \right] K_{j}u_{j},
\end{equation*}%
\begin{equation}
\left( O_{2j}-O_{2}\right) u_{j}=-\left( \frac{\partial \upsilon _{j}}{%
\partial x_{1}}-\frac{\partial \upsilon }{\partial x_{1}}\right) \frac{%
\partial }{\partial x_{1}}u_{j}+\left( \frac{\partial \upsilon _{j}}{%
\partial x_{2}}-\frac{\partial \upsilon }{\partial x_{2}}\right) \left(
\Lambda _{1}+\Lambda _{2}\right) u_{j},  \tag{7.6}
\end{equation}%
\begin{equation*}
\left( O_{3j}-O_{3}\right) u_{j}=\dsum\limits_{i=1}^{4}\left(
O_{3ij}-O_{3i}\right) u_{j}
\end{equation*}%
here 
\begin{equation*}
O_{3ij}=B_{0}\left( g\right) S\left( g\right) G_{i}\left( g\right) \left(
x_{1j},0\right) ,
\end{equation*}%
and $O_{3ij}$, $G_{i}$ are defined as in $\left( 6.14\right) $; moreover, 
\begin{equation*}
O\left( x_{1j}\right) =O\left( x_{1j},0\right)
=\dsum\limits_{k=1}^{3}O_{k}\left( g\right) \left( x_{1j}\right) ,\text{ }%
O_{k}\left( x_{1j}\right) =O_{k}\left( g\right) \left( x_{1j},0\right) ,%
\text{ }
\end{equation*}%
\begin{equation*}
B_{0j}=B_{0}\left( g\right) \left( x_{1j},0\right) ,\text{ }K_{j}=K\left(
g\right) \left( x_{1j},0\right) ,\text{ }\upsilon _{j}=\upsilon _{j}\left(
g\right) =K\left( g\right) \left( x_{1j},0\right) g,\text{ }
\end{equation*}%
\begin{equation*}
S_{j}=S\left( g\right) \partial B\left( g\right) \left[ .,K\left( g\right) g%
\right] \left( x_{1j},0\right) ,
\end{equation*}%
$O_{k}\left( x_{1j}\right) $ are local operators fixed at points $\left(
x_{1j},0\right) $\ defined by equalities $\left( 6.3\right) ,$ $\left(
6.5\right) ,$ $\left( 6.14\right) $, $\left( 6.15\right) $, respectively. \
From expressions of $F_{j},$ $K\left( g\right) \left( x_{1j},0\right) $ by
using $\left( 6.3\right) ,$ $\left( 6.5\right) ,$ $\left( 6.14\right) $, $%
\left( 6.15\right) $ we get that $F_{j}\in h^{1,\alpha }\left( A\right) .$

\bigskip For proving the Theorem 7.1 we need the following lemmas:

\textbf{Lemma 7.1. } The operator $u\rightarrow \left( O\left( x_{1j}\right)
+\mu ^{2}\right) u$ is an isomorphisim\ from $h^{2,\alpha }\left( A\right) $
onto $h^{1,\alpha }\left( A\right) .$ Moreover, the following estimate holds%
\begin{equation}
\dsum\limits_{i=0}^{1}\left\vert \mu \right\vert ^{2-i}\left\Vert \frac{%
\partial ^{i}u}{\partial x_{1}^{i}}\right\Vert _{h^{1,\alpha }\left(
A\right) }+\left\Vert A_{j}u\right\Vert _{h^{1,\alpha }\left( A\right) }\leq
C\left\Vert \left( O\left( x_{1j}\right) +\mu ^{2}\right) u\right\Vert
_{h^{1,\alpha }\left( A\right) }  \tag{7.7}
\end{equation}%
for all $u\in h^{2,\alpha }\left( A\right) $ and\ for sufficiently large $%
\mu >0$.

\textbf{Proof. }Consider the equation\textbf{\ }%
\begin{equation*}
\left( O\left( x_{1j}\right) +\mu ^{2}\right) u=f.
\end{equation*}%
Then, by virtue of Theorem 6.6 we obtain that the operator $u\rightarrow
\left( O\left( x_{1j}\right) +\mu ^{2}\right) u$ is an isomorphisim\ from $%
h^{2,\alpha }\left( A\right) $ onto $h^{1,\alpha }\left( A\right) $ and the
estimate $\left( 7.7\right) $ holds.

\bigskip \textbf{Lemma 7.2. }There is a positive $\varepsilon \in \left(
0,1\right) $ such that the following local estimate holds%
\begin{equation}
\left\Vert \left( O_{1j}-O_{1}\right) u_{j}\right\Vert _{h^{1,\alpha }\left(
A\right) }\leq \varepsilon \left\Vert K_{j}u_{j}\right\Vert _{h^{1,\alpha
}\left( A\right) }.  \tag{7.8}
\end{equation}%
\textbf{\ }

\bigskip \textbf{Proof. }\ From $\left( 7.6\right) $ we get 
\begin{equation*}
\left\Vert \left( O_{1j}-O_{1}\right) u_{j}\right\Vert _{h^{1,\alpha }\left(
A\right) }\leq \left\Vert B_{0}\left( g\right) \left[ K_{j}-K\left( g\right) %
\right] u_{j}\right\Vert _{h^{1,\alpha }\left( A\right) }+\left\Vert \left[
B_{0j}-B_{0}\left( g\right) \right] K_{j}u_{j}\right\Vert _{h^{1,\alpha
}\left( A\right) }.
\end{equation*}

Then by taking into account of expressions $B_{0}\left( g\right) $, $B_{0j}$%
, $K\left( g\right) $, $K_{j}$ and in view of smoothness of coefficients,
choosing $\delta $ sufficiently small\ we have 
\begin{equation*}
\left\Vert \left( O_{1j}-O_{1}\right) u_{j}\right\Vert _{h^{1,\alpha }\left(
A\right) }\leq \left\Vert B_{0}\left( g\right) \left[ K_{j}-K\left( g\right) %
\right] u_{j}\right\Vert _{h^{1,\alpha }\left( A\right) }+
\end{equation*}

\begin{equation*}
\left\Vert \left( b_{1}-b_{1j}\right) \frac{\partial }{\partial x_{1}}%
K_{j}u_{j}\right\Vert _{h^{1,\alpha }\left( A\right) }+\left\Vert \left(
b_{2}-b_{2j}\right) \frac{\partial }{\partial x_{2}}K_{j}u_{j}\right\Vert
_{h^{1,\alpha }\left( A\right) }\leq
\end{equation*}

\begin{equation*}
\varepsilon \left( \left\Vert \frac{\partial }{\partial x_{1}}%
K_{j}u_{j}\right\Vert _{h^{1,\alpha }\left( A\right) }+\left\Vert \frac{%
\partial }{\partial x_{2}}K_{j}u_{j}\right\Vert _{h^{1,\alpha }\left(
A\right) }\right) \leq \varepsilon \left\Vert K_{j}u_{j}\right\Vert
_{h^{2,\alpha }\left( A\right) }.
\end{equation*}

\bigskip \textbf{Lemma 7.3. }There is a positive $\varepsilon \in \left(
0,1\right) $ such that the following local estimate holds%
\begin{equation}
\left\Vert \left( O_{2j}-O_{2}\right) u_{j}\right\Vert _{h^{1,\alpha }\left(
A\right) }\leq \varepsilon \left\Vert u_{j}\right\Vert _{h^{1,\alpha }\left(
A\right) }.  \tag{7.9}
\end{equation}

\textbf{Proof. }\ From the $\left( 7.6\right) $ we get%
\begin{equation*}
\left\Vert \left( O_{2j}-O_{2}\right) u_{j}\right\Vert _{h^{1,\alpha }\left(
A\right) }\leq \left\Vert \left( \frac{\partial \upsilon _{j}}{\partial x_{1}%
}-\frac{\partial \upsilon }{\partial x_{1}}\right) \frac{\partial }{\partial
x_{1}}u_{j}\right\Vert _{h^{1,\alpha }\left( A\right) }+
\end{equation*}%
\begin{equation*}
\left\Vert \left( \frac{\partial \upsilon _{j}}{\partial x_{2}}-\frac{%
\partial \upsilon }{\partial x_{2}}\right) \left( \Lambda _{1}+\Lambda
_{2}\right) u_{j}\right\Vert _{h^{1,\alpha }\left( A\right) }.
\end{equation*}

Then, by using the smoothness of coefficients and choosing $\delta $
sufficiently small, we get the estimate $\left( 7.9\right) .$

\textbf{Lemma 7.3. }There is a positive $\varepsilon \in \left( 0,1\right) $
such that the following local estimate holds%
\begin{equation}
\left\Vert \left( O_{3j}-O_{3}\right) u_{j}\right\Vert _{h^{1,\alpha }\left(
A\right) }\leq \varepsilon \left\Vert u_{j}\right\Vert _{h^{1,\alpha }\left(
A\right) }.  \tag{7.10}
\end{equation}

\textbf{Proof. }From expressions $\left( 6.13\right) -\left( 6.14\right) $
we get 
\begin{equation}
\left\Vert \left( O_{3j}-O_{3}\right) u_{j}\right\Vert _{h^{1,\alpha }\left(
A\right) }\leq \dsum\limits_{i=1}^{4}\left\Vert \left( O_{3ij}-O_{3i}\right)
u_{j}\right\Vert _{h^{1,\alpha }\left( A\right) }.  \tag{7.11}
\end{equation}

Moreover, from expressions $O_{31}$ and $G_{1}$ in $\left( 6.14\right) $ by
boundedness of functions $g$, $\nu $, $\beta $ we have 
\begin{equation*}
\left\Vert \left( O_{31j}-O_{31}\right) u_{j}\right\Vert _{h^{1,\alpha
}\left( A\right) }\leq
\end{equation*}%
\begin{equation*}
C\left\Vert \left( B_{0}\left( g\right) S\left( g\right) \left(
x_{1j},0\right) -B_{0}\left( g\right) S\left( g\right) \left( x\right)
\right) \left( u_{j}+\frac{\partial u_{j}}{\partial x_{1}}\right)
\right\Vert _{h^{1,\alpha }\left( A\right) }\leq
\end{equation*}%
\begin{equation*}
C\left\{ \left\Vert \left[ B_{0}\left( g\right) S\left( g\right) \left(
x_{1j},0\right) -B_{0}\left( g\right) S\left( g\right) \left( x\right) %
\right] u_{j}\right\Vert _{h^{1,\alpha }\left( A\right) }\right. +
\end{equation*}%
\begin{equation*}
\left\Vert B_{0}\left( g\right) S\left( g\right) \left( x\right) \left[
u_{j}\left( x\right) -u_{j}\left( x_{1j},0\right) \right] \right\Vert
_{h^{1,\alpha }\left( A\right) }+
\end{equation*}%
\begin{equation*}
\left\Vert \left[ B_{0}\left( g\right) S\left( g\right) \left(
x_{1j},0\right) -B_{0}\left( g\right) S\left( g\right) \left( x\right) %
\right] \frac{\partial u_{j}}{\partial x_{1}}\right\Vert _{h^{1,\alpha
}\left( A\right) }+
\end{equation*}%
\begin{equation*}
\left\Vert B_{0}\left( g\right) S\left( g\right) \left( x\right) \left[ 
\frac{\partial u_{j}}{\partial x_{1}}-\frac{\partial }{\partial x_{1}}%
u_{j}\left( x_{1j},0\right) \right] \right\Vert _{h^{1,\alpha }\left(
A\right) }.
\end{equation*}

\bigskip Then by boundedness of operator $B_{0}\left( g\right) S\left(
g\right) ,$ smoothness of coefficients, choosing $\delta $ sufficiently we
obtain from the above 
\begin{equation}
\left\Vert \left( O_{31j}-O_{31}\right) u_{j}\right\Vert _{h^{1,\alpha
}\left( A\right) }\leq \varepsilon \left\Vert u_{j}\right\Vert _{h^{2,\alpha
}\left( A\right) }.  \tag{7.12}
\end{equation}

In a similar way, from expressions $O_{32}$ and $G_{2}$\ in $\left(
6.14\right) $ we have 
\begin{equation}
\left\Vert \left( O_{32j}-O_{32}\right) u_{j}\right\Vert _{h^{1,\alpha
}\left( A\right) }\leq \varepsilon \left\Vert u_{j}\right\Vert _{h^{2,\alpha
}\left( A\right) }.  \tag{7.13}
\end{equation}

In view of the condition on the operator $\partial A\left( g\right) ,$
boundedness of operator $B_{0}\left( g\right) S\left( g\right) ,$ smoothness
of coefficients, choosing $\delta $ sufficiently we obtain%
\begin{equation}
\left\Vert \left( O_{33j}-O_{33}\right) u_{j}\right\Vert _{D_{A}\left(
\alpha \right) }\leq \varepsilon \left\Vert u_{j}\right\Vert _{D_{A}\left(
\alpha \right) }.  \tag{7.14}
\end{equation}

Finally, from expressions $O_{34}$ an $G_{4}$ in $\left( 6.14\right) $ by
boundedness of functions $g$, $\nu $, $\beta $ we have%
\begin{equation*}
\left\Vert \left( O_{34j}-O_{34}\right) u_{j}\right\Vert _{h^{1,\alpha
}\left( A\right) }\leq
\end{equation*}%
\begin{equation*}
C\left\Vert \left( B_{0}\left( g\right) S\left( g\right) \left(
x_{1j},0\right) -B_{0}\left( g\right) S\left( g\right) \left( x\right)
\right) \left( u_{j}+\frac{\partial u_{j}}{\partial x_{1}}+\frac{\partial
^{2}u_{j}}{\partial x_{1}^{2}}\right) \right\Vert _{h^{1,\alpha }\left(
A\right) }\leq
\end{equation*}%
\begin{equation*}
C\left\{ \left\Vert \left[ B_{0}\left( g\right) S\left( g\right) \left(
x_{1j},0\right) -B_{0}\left( g\right) S\left( g\right) \left( x\right) %
\right] u_{j}\right\Vert _{h^{1,\alpha }\left( A\right) }\right. +
\end{equation*}%
\begin{equation}
\left\Vert B_{0}\left( g\right) S\left( g\right) \left( x\right) \left[
u_{j}\left( x\right) -u_{j}\left( x_{1j},0\right) \right] \right\Vert
_{h^{1,\alpha }\left( A\right) }+  \tag{7.15}
\end{equation}%
\begin{equation*}
\left\Vert \left[ B_{0}\left( g\right) S\left( g\right) \left(
x_{1j},0\right) -B_{0}\left( g\right) S\left( g\right) \left( x\right) %
\right] \frac{\partial u_{j}}{\partial x_{1}}\right\Vert _{h^{1,\alpha
}\left( A\right) }+
\end{equation*}%
\begin{equation*}
\left\Vert B_{0}\left( g\right) S\left( g\right) \left( x\right) \left[ 
\frac{\partial u_{j}}{\partial x_{1}}-\frac{\partial }{\partial x_{1}}%
u_{j}\left( x_{1j},0\right) \right] \right\Vert _{h^{1,\alpha }\left(
A\right) }+
\end{equation*}%
\begin{equation*}
\left\Vert \left[ B_{0}\left( g\right) S\left( g\right) \left(
x_{10},0\right) -B_{0}\left( g\right) S\left( g\right) \left( x\right) %
\right] \frac{\partial ^{2}u_{j}}{\partial x_{1}^{2}}\right\Vert
_{D_{A}\left( \alpha \right) }+
\end{equation*}%
\begin{equation*}
\left\Vert B_{0}\left( g\right) S\left( g\right) \left( x\right) \left[ 
\frac{\partial ^{2}u_{j}}{\partial x_{1}^{2}}-\frac{\partial ^{2}}{\partial
x_{1}^{2}}u_{j}\left( x_{1j},0\right) \right] \right\Vert _{B_{1p}}\leq
\varepsilon \left\Vert u_{j}\right\Vert _{h^{1,\alpha }\left( A\right) }.
\end{equation*}

Then the estiamate $\left( 7.10\right) $ is obtained from $\left(
7.12\right) -\left( 7.15\right) .$

Now, we can prove the Theorem 7.1.

\textbf{Proof of Theorem 7.1.} By virtue of Theorem 6.6 the following
estimate holds%
\begin{equation}
\dsum\limits_{i=0}^{1}\left\vert \mu \right\vert ^{2-i}\left\Vert \frac{%
\partial ^{i}u_{j}}{\partial x_{1}^{i}}\right\Vert _{h^{1,\alpha }\left(
A\right) }+\left\Vert A_{0}u_{j}\right\Vert _{h^{1,\alpha }\left( A\right)
}\leq C\left\Vert F_{j}\right\Vert _{h^{1,\alpha }\left( A\right) } 
\tag{7.16}
\end{equation}%
for all solution $u_{j}\in h^{2,\alpha }\left( A\right) $ of the equation $%
\left( 7.5\right) .$

Whence, using smoothness of coefficients of equations $\left( 7.3\right)
-\left( 7.6\right) ,$ in view of Lemmas 7.1-7.3 for $\mu $ with sufficiently
large Re$\mu >0$ we get%
\begin{equation}
\left\Vert F_{j}\right\Vert _{h^{1,\alpha }\left( A\right) }\leq \varepsilon 
\left[ \dsum\limits_{i=0}^{1}\left\vert \mu \right\vert ^{2-i}\left\Vert 
\frac{\partial ^{i}u_{j}}{\partial x_{1}^{i}}\right\Vert _{h^{1,\alpha
}\left( A\right) }+\left\Vert A_{j}u_{j}\right\Vert _{h^{1,\alpha }\left(
A\right) }\right] +\left\Vert f.\varphi _{j}\right\Vert _{h^{1,\alpha
}\left( A\right) }.  \tag{7.17}
\end{equation}

Moreover, by appling the microlocal analysis reasoning as in theorems 1, 2
in $\left[ 21\right] $ and in theorem 6.2 and as in $\left[ 9\text{, Lemma
6. 5-6.7}\right] $\ we obtain the same estimates for corresponding
commutators operators. Then from this and from $\left( 7.16\right) $, $%
\left( 7.17\right) $ we get the assertion.

From the Theorem 7.1 by microlocal analysis reasoning as in theorem 6.2 and
corollary 6.3 in $\left[ \text{9}\right] $ we obtain

\textbf{Corollary 7.1. }Assume the conditions 3.1, 3.2 are satisfied and $%
K\subset W$ is compact. Then there exist positive constants $\mu _{0}$ and $%
C=C(K)$ such that

\begin{equation*}
\left\Vert \upsilon \right\Vert _{h^{2,\alpha }\left( A\right) }+\left\vert
\mu \right\vert \left\Vert \upsilon \right\Vert _{h^{1,\alpha }\left(
A\right) }\leq C\left\Vert \left( \mu +\partial O_{t}\left( g\right) \right)
\upsilon \right\Vert _{h^{1,\alpha }\left( A\right) }
\end{equation*}%
for all $\upsilon \in h^{2,\alpha }\left( A\right) ,$ $g\in K,$ $t\in \left[
0,1\right] $ and $\mu \in \left\{ z:\func{Re}z>\mu _{0}\right\} .$

\textbf{Corollary 7.2. }Let $g\in W$ and $t\in \left[ 0,1\right] $ be given.
Then 
\begin{equation*}
\partial O_{t}\left( g\right) \in H\left( h^{2,\alpha },h^{1,\alpha }\right)
.
\end{equation*}

Now, by using Theorem 7.1 we can prove Theorem 1.

\textbf{Proof of Theorem 1: }Let $f_{0}\in V_{\nu }$ be given and set $%
g_{0}=f_{0}-\nu .$ Observe that $g_{0}\in W_{\alpha ,1}=W.$ It follows from
Lemma 4.1 that we only have to prove that there exist $t_{+}>0$ and a unique
maximal classical solution of $\left( 1.1\right) -\left( 1.2\right) $ on $%
[0,\left. t_{+}\right) $ satisfying 
\begin{equation}
\lim\limits_{t\rightarrow t^{+}}\left\Vert g\left( t,.\right) \right\Vert
_{h^{2,\alpha }}=\infty \text{, }\lim\limits_{t\rightarrow
t^{+}}\inf\limits_{\upsilon \in \partial W}\left\Vert f\left( t,.\right)
-\upsilon \right\Vert _{h^{2,\alpha }}=0  \tag{7.18}
\end{equation}%
if $t_{+}<\infty $ and $g\in C_{b}\left( [0,\left. t_{+}\right) ;W\right) .$

By reasoning as $\left[ \text{9, Lemma 5.10}\right] $ It follows that $W$ is
an open subset of $h_{+}^{2,\alpha }$. Hence, thanks to Lemma 4.2 and
Corollary 7.2, we know that $O\in H\left( W,h^{1,\alpha }\right) $ and that 
\begin{equation}
\partial O\left( g\right) \in H\left( h^{2,\alpha },h^{1,\alpha }\right) 
\text{, }g\in W.  \tag{7.19}
\end{equation}

\bigskip Let now $\beta \in \left( 0,\alpha \right) $ be fixed and observe
that $W\subset W_{\beta ,1}$. Thus the very same arguments as above also
ensure that%
\begin{equation*}
\partial O\left( g\right) \in H\left( h^{2,\beta },h^{1,\beta }\right) \text{%
, }g\in W.
\end{equation*}

It is not difficult to see that the maximal $h^{1,\alpha }$-realization of $%
\partial O\left( g\right) \in B\left( h^{2,\beta },h^{1,\beta }\right) $ for 
$g\in W$, is just the linear operator in $\left( 7.19\right) $. Note that 
\begin{equation*}
\left( h^{1,\beta },h^{2,\beta }\right) _{\alpha -\beta ,\infty
}=h^{1,\alpha }
\end{equation*}%
where $\left( .,.\right) _{\alpha -\beta ,\infty }$ denotes the real
interpolation. Consequently, invoking Theorem 2.3 in $\left[ 24\right] $, we
find that 
\begin{equation}
\partial O\left( g\right) \in M_{1}\left( h^{2,\alpha },h^{1,\alpha }\right) 
\text{, }g\in W,  \tag{7.20}
\end{equation}%
where $M\left( E_{1},E_{2}\right) $ denotes the class of all operators in s
having the property of maximal regularity in the sense of Da Prato and
Grisvard $\left[ 6\right] $. The existence of a unique maximal classical
solution of (E)g o and the property of a smooth semiflow on $W$ can now be
obtained along the lines of the proofs of Proposition 3.5 and Theorem 3.2 in 
$[24]$.

Finally suppose that $t_{+}<\infty $, $g\in C_{b}\left( [0,\left.
t_{+}\right) ;W\right) $ and that $\left( 7.18\right) $ is not true. Then $%
g_{1}=$ $\lim\limits_{t\rightarrow t_{+}}$ $g\left( t\right) $ exists in $W$%
. Hence taking \ \ \ $g_{1}$ as initial value in $\left( 4.2\right) $ one
easily constructs a solution $\bar{g}$\ of $\left( 4.2\right) $ for initial
date $g_{1}$ extending $g$. This contradicts the maximality of $g$.

\begin{center}
\ 
\end{center}

\ \textbf{References}

\begin{quote}
\ \ \ \ \ \ \ \ \ \ \ \ \ \ \ \ \ \ \ \ \ \ \ \ 
\end{quote}

\begin{enumerate}
\item Amann H., Maximal regularity for non-autonomous evolution equations.
Advanced Nonlinear Studies, (2004) 4, 417-430.\ \ \ \ \ \ \ \ \ \ \ \ \ \ \
\ \ \ \ \ \ \ \ \ \ \ \ \ \ \ \ \ \ \ \ \ \ \ \ \ \ \ \ \ \ \ \ \ \ \ \ \ \
\ \ \ \ \ \ \ \ \ \ \ \ \ 

\item Agarwal R., Bohner M., Shakhmurov V. B., Linear and nonlinear nonlocal
boundary value problems for differential operator equations, Appl. Anal.,
(2006), \ 85(6-7), 701-716.

\item Ashyralyev A, Cuevas. C and Piskarev S., On well-posedness of
difference schemes for abstract elliptic problems in spaces, Numer. Func.
Anal. Opt., (2008)29, (1-2), 43-65.\ 

\item Anosov V. P., Sobolevskii P. E., Coercive solvability boundary value
problems for elliptic equations of the second orders in Banach spaces, I,
II, Differ. Uravn., (7)11, 2030-2044, (7)12 (1971), 2191--2198.

\item Arendt W., Batty C., Bu S., Fourier multipliers for Holder continuous
functions and maximal regularity, Studia Mathematica 160 (1) (2004), 1-29.

\item Da Prato P., Grisvard, P., Equations d'\'{e}volution abstraites
nonlinair\'{e}s de type parabolique, Ann. Mat. Pura Appl., (4) 120, 329-396
(1979).

\item Dore, G., $L_{p}$-regularity for abstract differential equations. In:
Functional Analysis and Related Topics, H. Komatsu (ed.), Lecture Notes in
Math. 1540. Springer, 1993.

\item Denk R., Hieber M., Pr\"{u}ss J., $R$-boundedness, Fourier multipliers
and problems of elliptic and parabolic type, Mem. Amer. Math. Soc. (2003),
166 (788), 1-111.

\item Escher J., Simonett G., Maximal regularity for a free boundary
problem, NoDEA 2 (1995) 463-510.

\item Favini A., Shakhmurov V., Yakubov Y., Regular boundary value problems
for complete second order elliptic differential-operator equations in UMD
Banach spaces, Semigroup Form, (2009), 79 (1), 22-54.

\item Favini, A., Yagi, A., Degenerate Differential Equations in Banach
Spaces, Taylor \& Francis, Dekker, New-York, 1999.

\item Goldstain J. A., Semigroups of Linear Operators and Applications,
Oxford University Press, Oxfard, 1985.

\item M. Girardi and L. Weis, Operator-valued Fourier multiplier theorems on
Besov spaces, Mathematische Nachrichten, Math. Nachr. 251 (2003), 34-51.

\item Krein S. G., Linear Differential Equations in Banach space, American
Mathematical Society, Providence, 1971.

\item Lunardi A., Analytic Semigroups and Optimal Regularity in Parabolic
Problems, Birkhauser, 2003.

\item Lions J. L., Peetre J., Sur one classe d'espases d'interpolation, IHES
Publ. Math. (1964)19, 5-68.

\item Shklyar, A.Ya., Complete second order linear differential equations in
Hilbert spaces, Birkhauser Verlak, Basel, 1997.

\item Sobolevskii P. E., Coerciveness inequalities for abstract parabolic
equations, Dokl. Akad. Nauk, (1964), 57(1), 27-40.

\item Shahmurov R., On strong solutions of a Robin problem modeling heat
conduction in materials with corroded boundary, Nonlinear Anal. Real World
Appl., (2011),13(1), 441-451.

\item Shahmurov R., Solution of the Dirichlet and Neumann problems for a
modified Helmholtz equation in Besov spaces on an annuals, J. Differential
Equations, 2010, 249(3), 526-550.

\item Shakhmurov V. B., Maximal regular abstract elliptic equations and
applications, Siberian Mathematical Journal, V.51, no 5, 935-948, 2010.

\item Shakhmurov V. B., Degenerate differential operators with parameters,
Abstr. Appl. Anal., (2007), 2006, 1-27.

\item Shakhmurov V. B., Separable anisotropic differential operators and
applications, J. Math. Anal. Appl., 327(2) 2006, 1182-1201.

\item Simonett, G. Quasilinear parabolic equations and semiflows. In
Evolution Equations, Control Theory, and Biomathematics, Lecture Notes in
Pure and Appl. Math. M. Dekker, New York, 523-536, (1994).

\item Triebel H., Interpolation Theory, Function Spaces, Differential
Operators, North-Holland, Amsterdam, 1978.

\item Triebel H., Theory of Function Spaces II, Birkhuser-Verlag, Basel,
1992.

\item Weis L, Operator-valued Fourier multiplier theorems and maximal $L_{p}$
regularity, Math. Ann., (2001), 319, 735-758.

\item Yakubov S. and Yakubov Ya., Differential-Operator Equations. Ordinary
and Partial Differential Equations, Chapman and Hall /CRC, Boca Raton, 2000.
\end{enumerate}

\end{document}